\documentclass[12pt]{article}
\usepackage{amsfonts, amsmath,graphicx, amsthm, amssymb, cite, stmaryrd}
\numberwithin{equation}{section}
\newtheorem{proposition}{Proposition}
\numberwithin{proposition}{section}

\numberwithin{thm}{section}

\numberwithin{definition}{section}

\newtheorem{rem}{Remark}
\numberwithin{rem}{section}

\numberwithin{lem}{section}

\numberwithin{cor}{section}
\newcommand{\RR}{\mathrm{I\!R\!}}
\usepackage[margin=1in]{geometry}
\usepackage{color}

\usepackage{pgfplots}

\usetikzlibrary{arrows.meta}
\usetikzlibrary{shapes,backgrounds}
\usetikzlibrary{patterns}
\usetikzlibrary{decorations.markings}

\pgfplotsset{compat = newest}

\pgfplotsset{
  every axis/.append style={
    axis x line=middle,    % put the x axis in the middle
    axis y line=middle,    % put the y axis in the middle
    axis line style={<->,color=blue}, % arrows on the axis
    xlabel={$x$},          % default put x on x-axis
    ylabel={$y$},          % default put y on y-axis
  }
}

\begin{document}
%A semi-inverse problem in plane linear elastostatics.\\
\title{On singular behaviour in a plane linear elastostatics problem}
\author{H.~Gimperlein\thanks{Engineering Mathematics, University of Innsbruck, Innsbruck, Austria} \and M.~Grinfeld\thanks{Department of Mathematics and Statistics, University of Strathclyde, Glasgow, G1 1XH, UK} \and R.~J.~Knops\thanks{The Maxwell Institute of Mathematical Sciences and School of Mathematical and Computing Sciences, Heriot-Watt University, Edinburgh, EH14 4AS, Scotland, UK} \and M.~Slemrod\thanks{Department of Mathematics, University of Wisconsin, Madison, WI 53706, USA}}
\maketitle
\begin{center}
Dedicated to Professor Marcelo Epstein on the 80th anniversary of his birthday.
\end{center}
\begin{abstract}
\noindent  A vector field similar to  those separately introduced by Artstein and Dafermos is constructed from the tangent to a monotone increasing one-parameter family of non-concentric circles that touch at the common point of intersection taken as the origin. The circles define and space-fill a lens shaped region $\Omega$ whose outer and inner boundaries are the greatest and least circles. The double cusp at the origin creates a geometric singularity at which the vector field is indeterminate and has non-unique limiting behaviour. A semi-inverse method that involves the Airy stress function then shows that the vector field corresponds to the displacement vector field for   a linear plane compressible non-homogeneous isotropic elastostatic equilibrium problem  in $\Omega$ whose boundaries  are rigidly rotated relative to each other, possibly causing  rupture or tearing at the origin.  A sequence of   solutions is found for which not only are the Lam\'{e} parameters  strongly-elliptic, but the non-unique limiting behaviour of the displacement is preserved.  Other properties of the vector field are also established.  
\end{abstract}

\section{Introduction}\label{intro}

This paper continues the investigation  relevant to   continuum mechanics  
of a vector field $u(x,t)$    introduced separately  by Z. Artstein \cite{a83} and C. M. Dafermos \cite{d12}. The vector field,  defined on some bounded region $\Omega\subset\RR^{2}$,  has components   
\begin{eqnarray}
\label{u1}
u_{1}(x_{1},\,x_{2},t)&=& x_{2},\\
\label{u2}
u_{2}(x_{1},\,x_{2},t)&=& \frac{(x^{2}_{2}-x^{2}_{1})}{2x_{1}},
\end{eqnarray}
with respect to  Cartesian coordinates $(x_{1},\,x_{2})$  and time variable $t$. %A notable feature of \eqref{u2} is that it is  indeterminate at the coordinate origin $(0,\,0)$.

Artstein considers a variant of \eqref{u1} and \eqref{u2} appropriate to  a control system for which a smooth stabilizing control does not exist. Independently, Dafermos  used \eqref{u1} and \eqref{u2} with $u_{1}=\dot{x}_{1},\,u_{2}=\dot{x}_{2}$ as an example of an initial value problem for   an evolutionary ordinary differential equation that possesses a non-unique non-smooth solution. Uniqueness  is recovered on application of an ``entropy rate'' admissibility criterion. The present authors \cite{ggks23}  extend the study  of \eqref{u1} and \eqref{u2} by proving   that  the vector field $u(x,t)$  forms the Lagrangian trajectories of fluid particles  in steady  compressible flow. The trajectories   can be shown non-unique. The entropy rate criterion is again employed to  select a unique trajectory,  demonstrating the effectiveness of the criterion,

In this paper,  the vector field $u(x)$, supposed time independent, is derived in a manner different to \textcolor{black}{those of either Artstein or Dafermos.} As such, it is capable of obvious generalization. The vector field is taken parallel to the tangent at each point of a circle belonging to a monotone increasing one-parameter space-filling family whose members are contiguous at the coordinate origin. This leads immediately not only to the expressions \eqref{u1}, \eqref{u2} for the components $u_1$ and $u_2$ of $u$, but also as shown by \eqref{u2} to the indeterminacy of  $u_2$ at the origin. Moreover, the region of definition $\Omega$ is that enclosed by circles of the family of least and greatest radii which necessarily creates a double cusp at the origin. An important aim, therefore, of the paper is to explain this singular behaviour by seeking an interpretation in the specific context of a boundary value problem of plane linear elastostatics.

Before doing so, however, it is useful to establish further properties of the vector field $u$. Thus it is shown that the vector field has constant magnitude $c$ at each point of a member circle of radius $c$, that the limit of the component $u_2$ as the origin is approached along circles is ambiguous, and that the inner and outer boundary circles of $\Omega$ are rigidly rotated with respect to each other.

The last result provides the boundary conditions for the elastic problem, in which the vector field is assumed to be the prescribed displacement field. However, the elastic problem has still to be defined for which \eqref{u1} and \eqref{u2} are  the components of a displacement that corresponds to a stress distribution in equilibrium under zero body force. It is easily checked that the elastic body connot be homogeneous. When further restricted to be linear, isotropic and compressible, both the stress and material parameters are obtained by a semi-inverse method that involves the Airy stress function. The twisting of the inner and outer boundaries, and the singular behaviour at the origin, may then be interpreted as due to tearing, rupture, or brittle fracture. The conjecture is supported by calculating the force and couple acting on $\Omega$ together with the strain energy.

Section~\ref{not} introduces notation and reviews selected aspects of plane elastostatics. Section~\ref{props} derives the components of the  vector fields tangential to non-concentric circles belonging to a monotone increasing one-parameter family whose members all touch at the common point of intersection taken to be  the origin.  
Section \ref{base} describes the region of definition $\Omega$ and defines the vector field $u(x)$ of current interest as that parallel to the tangent vector field introduced in Section \ref{props}. Section \ref{base} also \textcolor{black}{explores further} properties of $u(x)$. In particular, at the origin, where the circles intersect and have a common tangent, the second component \eqref{u2} of the vector, besides being  indeterminate, exhibits   non-unique limiting behaviour. Section~\ref{asf} employs a semi-inverse method and Airy's stress function to prove that the vector field is the vector displacement field appropriate to  displacement boundary value problems on $\Omega$ for a nonhomogeneous compressible isotropic plane elastic body in equilibrium subject to zero body force and rigidly rotated boundaries. Conditions are derived for  the associated   Lam\'{e} parameters to be  strongly-elliptic. \textcolor{black}{Increased} understanding of the singular behaviour at the origin of the displacement is sought in Sections  \ref{fse} and \ref{gense}. Section \ref{fse} calculates that the total traction acting on $\Omega$ is zero whereas the total couple is proportional to the difference \textcolor{black}{between} the squares of the outer and inner circles. Section~\ref{gense} is devoted to calculating the corresponding total strain energy which under certain conditions does not have a singularity while in others the singularity varies inversely as the distance from the origin.  Section~\ref{concl} presents certain directions in which the analysis might be extended and lists some open problems. An Appendix derives the Airy stress function used in Section~\ref{asf} and proposes a generalisation.\\

Vector and scalar quantities are not typograhically distinguished although subscripts denote the components of a vector. The subscript comma and summation conventions are adopted throughout with Latin suffixes ranging over $1,2,3$ while Greek suffixes have the range $1,2$. Since the vector field $u$ is prescribed it immediately follows that  the solution to the elastic boundary value problem always exists.

\section{Notation, basic theory  and other preliminaries}\label{not}

Several properties of the vector field leading to \eqref{u1} and \eqref{u2} do not depend upon the  relation to elasticity or even continuum mechanics. Later sections, however, require the vector field to be interpreted as the displacement vector occurring in linear plane elasticity. Consequently, this Section, apart from notation, introduces notions from elasticity theory in preparation for the subsequent discussion. The reader is referred to e.g., \cite{g72,ll70,l27} and other standard texts for complete accounts.

Let $\Omega\subset \RR^{2}$   denote a bounded simply or multiply connected  plane region and  let $(x_{1},\,x_{2})$ be the Cartesian  coordinates of a point $x\in\Omega$. Suppose $\Omega$ is occupied by  material that   deforms to the region   $\widehat{\Omega}$ and that the particle  at $x$ is displaced to  the position $y\in \widehat{\Omega}$. The map $\Omega \rightarrow \widehat{\Omega}$, represented  by $y=y(x)$, is supposed time independent and sufficiently differentiable. The displacement vector field  $u(x)$ (c.p., \cite{ll70,g72}) is defined as
%%and that the point $x\in\Omega$ is deformed to the point $y\in\widehat{\Omega}$ in the deformed body $\widehat{\Omega}$.  . The vector field $u(x)$   
%%Let $u(x)$ be a time-independent vector field whose Cartesian components are \eqref{u1} and \eqref{u2}. Features of $u$ regarded as a  displacement vector in continuum mechanics are reviewed in Section~\ref{el}, but our immediate task is to explain the choice \eqref{u1} and \eqref{u2}. Although Artstein \cite{\a83} and Dafermos \cite{d12} have separately motivated this  choice, a different approach is adopted here based upon a one parameter family of closed simple curves.
%%This Section defines the basic domain of definition required in the discussion of  general properties of $u$ valid irrespective of whether of not $u$ is the displacement field of a particular elastic boundary value problem. However some notions for this problem are introduced preparatory to a fuller account in Section#\ref{el}.  
% %whose boundary where it exists may be arbritarily rough until otherwise stated.  %and let $x\in \Omega$ be some (mathematical) point(particle)  in $\Omega$ that is not to be confused with the fundamental particles of the body subsequently  modelled.% Due to various internal and external factors, to be specified, 
%%When regarded as a  displacement,
\[
u(x)=y(x)-x,
\]
so that  $u(x)$   as a function of $x$ completely determines $y(x)$.
%%Denote by  $\nabla u(x)$  the  displacement gradient of $u$  and
Recall that
\[
du_{i}=u_{i,j}dx_{j}.
\]
% and that the square of the deformed element of distance is given by
%\begin{eqnarray}
%\nonumber
%dy_{i}dy_{i}&=& (du_{i}+dx_{i})(du_{i}+dx_{i})\\
%\nonumber
%&=& dx_{i}dx_{i}+2u_{i,j}dx_{i}dx_{j}+u_{i,k}u_{i,j}dx_{k}dx_{j}\\
%\label{extel}
%&=& dx_{i}dx_{i}+2\eta_{ij}dx_{i}dx_{j},
%\end{eqnarray}
%where the rectangular Cartesian symmetric components $\eta_{ij}$  of  the nonlinear strain tensor are given by
%\begin{equation}
%\label{strdef}
%\eta_{ij}=(1/2)\left(u_{i,j}+u_{j,i}+u_{p,i}u_{p,j}\right).
%\end{equation}

Linear theories assume that the  displacement gradient is small 
and that to first order the Cartesian components of 
the strain are expressed  by
\begin{equation}
\label{linstradef}
e_{ij}=(1/2)\left(u_{i,j}+u_{j,i}\right).
\end{equation}
While the  vector field $u(x)$ is not necessarily correspondingly small,  the change $\delta V$  in the infinitesimal volume element at a given point may  be shown to be
\[
\delta V=e_{ii},
\]
which  may be either positive (extension), negative (compression) or zero (incompressible).

%Attention is now confined to linear compressible isotropic plane elastostatics for which %%%%%%%%%%%%%%%%%%%%%Here030624
%%%%%5\section{Plane elastostatics}\label{el}
%We briefly sketch for convenience the basic theory  of linear elastostatics that holds for displacements admitting singularities. Point defects and similar singular behaviour require the standard theory to be modified since in their neighbourhood the strain tensor becomes unbounded and the linear theory may be  no longer valid.  %Furthermore, use of  the divergence theorem requires  careful consideration.  

For a compressible linear isotropic plane elastic body, the symmetric stress components satisfy  the constitutive relations
\begin{equation}
\label{ss}
\sigma_{\alpha\beta}=\lambda e_{\gamma\gamma}\delta_{\alpha\beta}+2\mu e_{\alpha\beta},
\end{equation}
where $\alpha,\,\beta=1,\,2$,  the Kronecker delta function is denoted by $\delta_{\alpha\beta}$, and $\lambda(x),\,\mu(x)$ are the Lam\'{e} parameters related to Poisson's ratio $\nu(x)$ by 
\begin{equation}
\label{lapo}
\lambda= \frac{2\mu\nu}{(1-2\nu)},\qquad \nu\neq 1/2,
\end{equation}
or  
\begin{equation}
\label{pr}
2\nu= \frac{\lambda}{(\lambda+\mu)},\qquad \nu\neq 1/2.
\end{equation}

%For an incompressible elastic body, the strain has zero dilatation so that
%\begin{equation}
%\label{incomstr}
%e_{\alpha\alpha}=u_{\alpha,\alpha}=0,
%\end{equation}
%and Poisson's ratio is given by $\nu=1/2$. The constitutive relations \eqref{ss}  are  modified to 
%\begin{equation}
%\label{incomss}
%\sigma_{\alpha\beta}=p\delta_{\alpha\beta}+2\mu e_{\alpha\beta},
%\end{equation}
%where $p(x,y)$ is the unknown pressure.
%%%%%%%%%%%%%%%%%%%%%%%

The following relations are easily deduced from \eqref{ss}:
\begin{eqnarray}
\label{murel}
2\mu&=& \frac{\sigma_{12}}{e_{12}}=\frac{(\sigma_{11}-\sigma_{22})}{(e_{11}-e_{22})},\\
\label{lamrel}
2(\lambda+\mu)&=& \frac{(\sigma_{11}+\sigma_{22})}{(e_{11}+e_{22})},\\
\label{lam2mu}
(\lambda+2\mu)&=& \frac{(\sigma_{11}e_{11}-\sigma_{22}e_{22})}{(e^{2}_{11}-e^{2}_{22})},\\
\label{lamonl}
\lambda&=& \frac{(\sigma_{22}e_{11}-\sigma_{11}e_{22})}{(e^{2}_{11}-e^{2}_{22})}.
\end{eqnarray} 
%%%%%%%%%%%%%%%%%%%%%%%%%%%%%%%%%%%%%5
From  relations \eqref{murel}, we conclude  that the parameter $\Lambda$, defined by
\begin{equation}
\label{Lamdefone}
\Lambda (x_{1},\,x_{2}):= \frac{(e_{11}-e_{22})}{e_{12}}=\frac{(\sigma_{11}-\sigma_{22})}{\sigma_{12}},
\end{equation}
is explicitly independent of the Lam\'{e} parameters in both  nonhomogeneous and homogeneous elastic bodies.
%%%%%%%%%%%%%%%%%%%%%%220524one%%%%%%%%%%%%%%%%

The stress tensor  in equilibrium under zero body force has zero divergence, and therefore satisfies  
\begin{equation}
\label{seq}
\sigma_{\alpha\beta,\beta}=0.
\end{equation}
A general solution to these equilibrium equations, irrespective of constitutive relations,  may be expressed in terms of the Airy stress function $\chi(x_{1},\,x_{2})$ as
\begin{equation}
\label{sfn}
\sigma_{11}=-\chi_{,22},\qquad \sigma_{12}=\chi_{,12},\qquad \sigma_{22}=-\chi_{,11}.
\end{equation}

For  compressible non-homogeneous linear elasticity, substitution of \eqref{ss} in  the equilibrium equations \eqref{seq} establishes that the vector field $u(x)$ at place $x$ satisfies the Navier equations
\begin{equation}
\label{innav}
\mu u_{\alpha,\beta\beta}+(\lambda+\mu) u_{\beta,\beta\alpha}+\lambda_{,\alpha}u_{\beta,\beta}+\mu_{,\beta}(u_{\alpha,\beta}+u_{\beta,\alpha})=0,\qquad x\in\Omega,
\end{equation}
which for homogeneous elasticity becomes
\begin{equation}
\label{nav}
\mu u_{\alpha,\beta\beta}+(\lambda+\mu) u_{\beta,\beta\alpha}=0,\qquad x\in\Omega.
\end{equation}

%The corresponding systems for an incompressible linear isotropic elastic body are
%\begin{eqnarray}
%\label{inhomincomnav}
%\mu u_{\alpha,\beta\beta}+\mu_{,\beta}(u_{\alpha,\beta}+u_{\beta,\alpha})+p_{,\alpha}&=&0, \qquad x\in\Omega,\\
%\label{incomnav}
%\mu u_{\alpha,\beta\beta}+p_{,\alpha}&=&0,\qquad x\in\Omega.
%\end{eqnarray}

For a bounded domain $\Omega$ with piecewise smooth boundary $\partial\Omega$ and subject to standard Dirichlet boundary conditions,   \eqref{innav} admits a unique classical or weak solution if the strong-ellipticity condition holds (cf., \cite[pp. 19 and 62]{kp71});
\begin{equation}
\label{uniqq}
\mu(\lambda+2\mu)>0.
\end{equation}
This may be equivalently written  as 
\begin{equation}
\label{pruniq}
-\infty< \nu < 1/2,\qquad 1<\nu <+\infty,\qquad \mu\neq 0.
\end{equation}

\section{Unit tangent vector field}\label{props}

The vector field $u(x)$ whose components are  \eqref{u1} and \eqref{u2} is obtained by  assuming it to be  parallel to a certain unit tangent field.  The  argument,  motivated by  \cite[App. D]{ggks23}, is  different from that of either Artstein \cite{a83} or Dafermos \cite{d12}.%  motivates the introduction of the vector field \eqref{u1} and \eqref{u2}. %The argument  is generalised  in the concluding section. 

Introduce  a  one parameter family of non-concentric circles centred at $(c,0)$, of radius $c$, and  denoted by  
%%%%%%%%%%%%%%%%ok220424
\begin{equation}
 \label{phidef}
\psi_{c}(x):=c^{-2}\left\{(x_{1}-c)^{2}+x_{2}^2\right\}=1
\end{equation}
where the constant $c$ satisfies $1\le c\le R$ for constant $R$. The circles touch at the origin $(0,0)$ which is their common single point of intersection.  %More generally, we consider the family of closed smooth simple curves  
%\begin{equation}
%\label{famdef}
%\psi_{c}(x)=1,
%\end{equation}
%assumed likewise contiguous at the origin.
 
%%%%%%%%%%%%%%%%%%50o220424

The unit tangent vector  at the point $x$ on the circle \eqref{phidef}  has Cartesian components
\begin{eqnarray}
\label{t1}
t_{1}(x_{1},\,x_{2})&=& \frac{c}{2}\frac{\partial\psi_{c}}{\partial x_{2}}=x_{2}/c,\\
\label{t2}
t_{2}(x_{1},\,x_{2})&=& -\frac{c}{2}\frac{\partial\psi_{c}}{\partial x_{1}}=(c-x_{1})/c,
\end{eqnarray}
so that at the origin all curves have the common unit tangent whose components are given by  $(t_{1}(0,0),\,t_{2}(0,0)=(0,1).$  Note that for each point  on the circle $\psi_{c}=1$  \emph{except}  the origin, appeal to   \eqref{phidef} enables \eqref{t2} to be alternatively written as
\begin{eqnarray}
\nonumber
t_{2}(x_{1},\,x_{2})&=& \frac{2x_{1}(c-x_{1})}{2cx_{1}},\qquad x_{1}\ne 0,\\
%\nonumber
%&=& \frac{(2cx_{1}-2x^{2}_{1})}{2cx_{1}}\\
%\nonumber
%&=& \frac{\left(-(x_{1}-c)^{2}-x^{2}_{1}+c^{2}\right)}{2cx_{1}}\\
\label{t22def}
&=&\frac{(x^{2}_{2}-x^{2}_{1})}{2cx_{1}}.
\end{eqnarray} 
%%%%%%%%%%%%%%%%%%%%%555ok220424 

In the next section, the  vector field  $u$ is derived from the tangent components \eqref{t1} and \eqref{t2}.  Certain properties are established irrespective of whether $u(x)$ is interpreted as a displacement vector field. 
%%%%%%%%%%%%%%55ok290524

\section{Properties of the vector field}\label{base}
 %Before computing  tangent vectors to the circles   \eqref{phidef}, 
We now specify the shape of   the bounded region  $\Omega$.  Let the   outer and inner boundaries $\partial \Omega_{\alpha},\, \alpha=1,2$ of $\Omega$ be  determined by $\psi_{R}(x)=1$ and $\psi_{1}(x)=1$, where $R>1$.  Consequently, the outer boundary is the  circle centred at $(R,0)$ and of radius $R$ while the inner boundary is the  circle centred at $(1,0)$ and of radius $1$. The region  $\Omega$ is  lens shaped  with  symmetrically placed cusps at the origin which  therefore create geometric point singularities. Observe that the circles belonging to  \eqref{phidef}    fill  the whole  space occupied by   $\Omega$.  Each point $x\in\Omega$ belongs to one and only one member of the family \eqref{phidef}.

The vector field  $u$ is taken parallel to the  unit tangent vector $t$. % In particular, let the constants  $c$  satisfy $1\le c \le R$ and %suppose $w(x)$ is some scalar function. let $\psi_{c}$ be defined by \eqref{phidef}.
Set $u=ct$ to obtain at the point $(x_{1},\,x_{2})$ of the curve \eqref{phidef},  the expressions \eqref{u1} and \eqref{u2} which as shown in deriving \eqref{t22def}  may be written,  except at the origin, in the form 
%%%%%%%%%%%%%%%%%%%ok290524
\begin{eqnarray}
\label{u111}
u_{1}(x_{1},\,x_{2})&=& x_{2},\\
\label{u221}
u_{2}(x_{1},\,x_{2})&=& \frac{(x^{2}_{2}-x^{2}_{1})}{2x_{1}}\\
\label{u222}
&=& (c-x_{1}).
\end{eqnarray}
%Observe that  expressions \eqref{u111} and \eqref{u221} are   valid at each point   $x\in\Omega$ and not only for those $x$ restricted to describe a given member of the family \eqref{phidef}.
%%%%%%%%%%%%%%%%%%ok290524
%%%%%%%%%%%%%%%%%030624here
  
Components of the gradient when the vector field has components \eqref{u111} and \eqref{u221} are given by 
\begin{eqnarray}
\label{u1grad}
u_{1,1}&=& 0,\qquad u_{1,2} =1,\\
\label{u2grad}
u_{2,1}&=& -\frac{(x^{2}_{1}+x^{2}_{2})}{2x^{2}_{1}},\\
\label{u22grad} 
  u_{2,2}&=&\frac{x_{2}}{x_{1}},
\end{eqnarray}
 while the symmetric components, which by construction are compatible,  become
%\begin{eqn
%\label{eta11}
%\eta_{11}&=& \frac{(x^{2}_{1}+x_{2}^{2})^{2}}{4x^{4}_{1}}= \frac{1}{8\cos^{4}{\theta}},\\
%\label{eta12}
%\eta_{12}&=&\frac{(x_{1}^{2}-x_{2}^{2})}{4x_{1}^{2}}-\frac{x_{2}(x_{1}^{2}+x_{2}^{2})}{4x_%{1}^{3}}\\ 
%\nonumber
%&=&\frac{1}{4x^{3}_{1}}\left(x^{3}_{1}-x^{2}_{1}x_{2}-x_{1}x^{2}_{2}-x^{3}_{2}\right),\\
%\nonumber
%&=&\frac{1}{4}\left(1-\tan{\theta}-\tan^{2}{\theta}-\tan^{3}{\theta}\right)\\
%\label{eta22}
%\eta_{22}&=& \frac{1}{2x^{2}_{1}}\left(x_{1}+x_{2}\right)^{2}\\
%\nonumber
%&=& \frac{1}{2}\left(1+\tan{\theta}\right)^{2},
%\end{eqnarray}
%where polar coordinates $(r,\theta)$ with origin at $(0,0)$ satisfy the relations
%\[
%x_{1}=r\cos{\theta},\qquad x_{2}=r\sin{\theta}.
%The  linear strain components are   
\begin{eqnarray}
\label{e11}
e_{11}&=&0,\qquad e_{22}=\frac{x_{2}}{x_{1}},\\
\label{e12}
e_{12}&=& \frac{1}{4x_{1}^{2}}\left(x_{1}^{2}-x^{2}_{2}\right).
%\nonumber
%&=& \frac{1}{4}\left(1-\tan^{2}{\theta}\right).
\end{eqnarray}

%It also is convenient to express   the two-dimensional linear strain components   in polar coordinates $(r,\,\theta)$. We have:
%\begin{eqnarray}
%\label{polgrad1}
%e_{rr}=u_{r,r},& &\qquad e_{\theta\theta}=r^{-1}\left(u_{\theta,\theta}+u_{r}\right),\\
%\label{shgrad2}
%e_{r\theta}&=& (1/2)\left(u_{\theta,r}-r^{-1}u_{\theta}+r^{-1}u_{r,\theta}\right).
%\end{eqnarray} 
%%%%%%%%%%%%%%%%%%%%%%%%ok220524
%%%%%%%%%%%%%%%%%%%%%%%040624herehelm

%On selecting $w(x)=1$ we recover \eqref{u1} and \eqref{u2}. 

%\begin{rem}\label{w}
%An arbitrary scalar function $w(x)$ may be retained but for simplicity we take $w(x)=1$ in what follows. The general displacement $wu$,however,  is discussed in Section~\ref{concl} especially with respect to  point forces and Volterra dislocations in  homogeneous isotropic elastostatics.
%\end{rem}
%%%%%%%%%%%%%%%%%%%%%%%%%%ok290524

%\subsection{Continuous displacement}\label{cont}
%An interpretation of the vector field $u$ is now obtained that does not  rely upon  differentiability. This aspect is considered later. 

%%%%%%%%%%%%%%%%%%%ok220424

For given $c$ select the origin of the polar coordinate system $(\rho,\phi),\, -\pi \le \phi\le  \pi$ to be at the centre of the circle $\psi_{c}=1$.   A second polar coordinate  system $(r,\theta),\, -\pi/2\le \theta\le \pi/2$ based on the  origin  $(0,0)$  of the Cartesian system is related to the first by  
\begin{eqnarray}
\label{x1r}
x_{1}&=&r\cos{\theta}=c+\rho\cos{\phi},\\
\label{x2r}
 x_{2}&=&r\sin{\theta}=\rho\sin{\phi}.
\end{eqnarray}
Accordingly, the circle  \eqref{phidef} may be written as
\begin{equation}
\label{fampol2}
r=2c\cos{\theta},
\end{equation}
while the polar coordinates $(\rho,\phi)$ of points on the circle  satisfy 
\begin{equation}
\label{phitheta}
\rho=c,\qquad \phi=2\theta,
\end{equation}
and we have the relations 
\begin{eqnarray}
\label{famx}
x_{1}&=& r\cos{\theta}=2c \cos^{2}{\theta}=c \cos{2\theta}+c=c(1+\cos{\phi}),\\
\label{famy}
x_{2}&=& r\sin{\theta}=2c\cos{\theta}\sin{\theta}=c\sin{2\theta}=c\sin{\phi}.   \end{eqnarray}                                                                                             

In an obvious notation, the radial and tangential components of the vector field whose Cartesian components are \eqref{u111} and \eqref{u222} on using \eqref{fampol2} and \eqref{phitheta} in terms of polar coordinates $(\rho,\phi)$ are given by
\begin{eqnarray}
\nonumber
u_{\rho}(c,\phi)&=& u_{1}(x_{1},x_{2})\cos{\phi}+u_{2}(x_{1},x_{2})\sin{\phi}\\
\nonumber
&=& \left[r\sin{\theta}\cos{\phi}+(c-r\cos{\theta})\sin{\phi}\right]\\
\nonumber
&=& \left[-r\sin{(\phi-\theta)}+c\sin{\phi}\right]\\
\label{rad}
&=& 0,
\end{eqnarray}
and
\begin{eqnarray}
\nonumber
u_{\phi}(c,\phi)&=& -u_{1}(x_{1},x_{2})\sin{\phi}+u_{2}(x_{1},x_{2})\cos{\phi}\\
\nonumber
&=& \left[-r\sin{\theta}\sin{\phi}+(c-r\cos{\theta})\cos{\phi}\right]\\
%&=& \left[-r\cos{(\phi-\theta)}+c\cos{\phi}\right]\\
%\nonumber
%&=& \left[-r\cos{\theta}+c\cos{2\theta}\right]\\
\nonumber
&=& \left[-2c\cos^{2}{\theta}+c(2\cos^{2}{\theta}-1)\right]\\
\label{tan}
&=& -c.
\end{eqnarray}
%%%%%%%%%%%%%%%%%%%%%%%%%%ok220424
We conclude that the vector $u$ is of constant magnitude $c$ and is  tangential to the circle in the clockwise direction.

In a mechanical interpretation the boundary circles, for which $c=1$ and $c=R$,  rotate uniformly clockwise by different amounts $1, R$. The singular behaviour in the component \eqref{u221} at the origin, where all members of the family \eqref{phidef} touch, perhaps  corresponds to rupture or brittle fracture.  The infinitesimal volume element  in the neighbourhood of the origin is in tension when $x_{2}>0$ and in compression when $x_{2}<0$.

The relative rotation of the boundary circles  may be alternatively derived directly from \eqref{u111} and \eqref{u221} on expressing the vector field as 
\[
u(x_{1},x_{2})=u_{1}e_{1}+u_{2}e_{2},
\]
where $(e_{1},\,e_{2})$ are the base vectors in the Cartesian coordinate system whose origin is $(0,0)$. Substitution from \eqref{u111} and \eqref{u222} leads to 
\begin{equation}
\label{rigid}
u(x_{1},x_{2})= ce_{2}+(x_{2}e_{1}-x_{1}e_{2}),
\end{equation}
which shows that the vector field at each point on a circle  belonging to 
 \eqref{phidef} may be interpreted as  a rigid body displacement consisting of a rigid body translation of amount $c$ in the $x_{2}$-direction plus a rigid body unit rotation about the origin in a clockwise direction.

At the origin the vector field with components \eqref{u111} and \eqref{u222} is directed along the common tangent to the family of circles \eqref{phidef}  i.e., in the direction parallel to $(0,1)$ but the  magnitude is given by the parameter $c$ and thus depends upon the particular circle along which the origin is approached. In this sense, the vector field is not  uniquely defined at the origin. Indeed, we have proved the  following result:
\begin{proposition}[Limiting behaviour]\label{u2nonuniq}
For points on a given circle \eqref{phidef}, the limit of the component $u_{2}(x)$ specified by \eqref{u222} as the point $(x_{1},\,x_{2})$ tends to the origin depends along which member of the family \eqref{phidef} the origin is approached; that is   
\begin{eqnarray}
\label{limu22}
\lim_{x_{1}\rightarrow 0}{u_{2}(x_{1},x_{2})}=c,
\end{eqnarray}
where the constant $c$ can be arbitrarily chosen in the range  $1\le c \le R$.
Moreover, the  component $u_{2}(x)$  when given by \eqref{u221} is indeterminate at the origin   and  is not the limit of $u_{2}(x)$ defined on   any sequence of points belonging to a circle of the family \eqref{phidef}.
\end{proposition}

%The vector displacement field therefore does not have a  unique 
% limiting value as the origin is approached.% from within $\Omega$. 
The origin may be regarded as the limit of points belonging, say, to the different circles $\psi_{c_{1}}=1$ and $\psi_{c_{2}}=1$. We determine the jump  in the limit of $u_{2}(x)$. Consider  a point belonging to $\psi_{c_{2}}=1$ in the neighbourhood of the origin. Starting from this point and integrating around the separate circle $\psi_{c_{1}}=1$ yields 
\begin{eqnarray}
  \nonumber
  \left[u_{2}\right]&=& \lim_{x_{2}\rightarrow 0-}{u_{2}(x_{1},x_{2})}-\lim_{x_{2}\rightarrow 0+}{u_{2}(x_{1},x_{2})}\\
  \label{jump}
&=& c_{1}-c_{2}.
\end{eqnarray}
We conclude that the jump in  $u_{2}(x)$ in the limit as the origin is approached is ambiguous.

In general, differentiation of the vector field must be performed on the components \eqref{u111} and \eqref{u222} with the condition that the point $(x_{1},\,x_{2})$
belongs to a particular member  of the family \eqref{phidef} applied after  differentiation. 

Notwithstanding the last remark, observe that  the  vector  field  $u$ with components given by \eqref{u111} and either \eqref{u221} or \eqref{u222} is continuously differentiable at each point  in $\Omega$ except at the origin where,  as previously stated, it  becomes indeterminate.  By construction, at each point $x\in \Omega\backslash \{(0,\,0)\}$ the vector field is directed along the tangent to the appropriate circle $\psi_{c}=1$ and is of magnitude $c$. Intuitively, the vector field at  all points on a given member of \eqref{phidef} is \textcolor{black}{isochoric}, when interpreted as a displacement.

The result is proved on converting to the polar coordinate sytem $(r,\theta)$ introduced in \eqref{x1r} and \eqref{x2r}  for which  the radial and transverse components to \eqref{phidef} of the vector field are
\begin{eqnarray}
\nonumber
u_{r}&=& u_{1}\cos{\theta}+u_{2}\sin{\theta}\\
\nonumber
&=& r\sin{\theta}\cos{\theta}+(c-r\cos{\theta}) \sin{\theta}\\
\label{rad}
&=& c\sin{\theta},\\
\nonumber
u_{\theta}&=& -u_{1}\sin{\theta}+u_{2}\cos{\theta}\\
\nonumber
&=& -r\sin^{2}{\theta}+(c-r\cos{\theta})\cos{\theta}\\
\label{transs}
&=& -r+c\cos{\theta}.
%\label{transs}
%&=& -c\cos{\theta}.
\end{eqnarray} 
%%%%%%%%%%%%%%%% %%%%%%%%%%%%%%%%%%ok220424
The corresponding symmetric components of the gradient are given by
\begin{eqnarray}
  \nonumber
  e_{rr}=u_{r,r},\qquad & & \qquad e_{\theta\theta}= r^{-1}\left(u_{\theta,\theta}+u_{r}\right),\\
  \nonumber
  e_{r\theta}&=& \left(1/2\right)\left(u_{\theta,r}-r^{-1}u_{\theta}+r^{-1}u_{r,\theta}\right).
  \end{eqnarray}

Insertion of \eqref{rad} and \eqref{transs}  leads to  the infinitesimal dilatation in  the direction of a given circle  \eqref{phidef}   becoming 
\begin{eqnarray}
\nonumber
e_{rr}+e_{\theta\theta}&=& \frac{\partial u_{r}}{\partial r}+\frac{u_{r}}{r}+\frac{1}{r}\frac{\partial u_{\theta}}{\partial\theta}\\
\nonumber
&=& 0+\frac{c\sin{\theta}}{r}+\frac{1}{r}\left(-c\sin{\theta}\right)\\
\nonumber
&=& 0.
\end{eqnarray}                                                                  
Therefore, the vector field when its  gradient is  infinitesimal  may be interpreted as \textcolor{black}{isochoric} along each circle  $\psi_{c}(x_{1},x_{2})=1$ of the space-filling family.                                                                                    %Properties of the vector field \eqref{u111} and \eqref{u222} become of increased  significance  when it can be established that the vector field is the solution to a mathematical model of some physical process
The next Section  establishes that the vector field $u(x)$ is a displacement field   in  a linear plane elastostatics problem. %As stated in Remark~\ref{w}, homogeneous isotropic elastic bodies require special choices of $w(x)$ other than $w=1$, or  the introduction of appropriate body forces to ensure satisfaction of the Navier equation. Appropriate  choices of $w(x)$, however, may lead to  relationships with point defects such as Volterra and other dislocations. % which may be due to thermoelastic or electromagnetic factors.
%  This aspect is not discussed here.
  
%%%%%%%%%%%%%%%%%ok220524 

\section{Airy stress function}\label{asf}
We now prove that the vector field with Cartesian componenets \eqref{u1} and \eqref{u2} is the vector displacement that  occurs in a linear non-homogeneous isotropic compressible plane elastic body in equilibrium under zero  body-forces. The body must be non-homogeneous otherwise the Navier equations \eqref{nav} are not satisfied without addition of suitable  body forces. Nevertheless, subject to our assumptions, the stress has zero divergence and consequently may be expressed  in terms of Airy's stress function $\chi(x_{1},\,x_{2})$. Substitution of \eqref{sfn} in \eqref{Lamdefone} yields
\begin{equation}
\label{airy}
\chi_{,11}-\chi_{,22}-\Lambda(x_{1},\,x_{2})\chi_{,12}=0,\qquad \Lambda:=\frac{(e_{11}-e_{22})}{e_{12}},
\end{equation}
which represents the partial differential equation for $\chi$ that is explicitly independent of the Lam\'{e} parameters.

Although equation \eqref{airy} is hyperbolic and therefore amenable to  the method of characteristics, a solution may also be derived on transforming to the complex variable $z=x_{1}+ix_{2}$ and its conjugate $\bar{z}=x_{1}-ix_{2}$ where $i=\sqrt{-1}$. The respective partial derivatives are defined as
\begin{eqnarray}
\label{zdif}
\frac{\partial}{\partial z}&:=&\frac{1}{2}\left(\frac{\partial}{\partial x_{1}}-i\frac{\partial}{\partial x_{2}}\right),\\
\label{zbardif}
\frac{\partial}{\partial \bar{z}}&:=&\frac{1}{2}\left(\frac{\partial}{\partial x_{1}}+i\frac{\partial}{\partial x_{2}}\right),
\end{eqnarray}
and  the linear strain components \eqref{e11} and \eqref{e12} are
\begin{eqnarray}
\label{come11}
e_{11}(z,\bar{z})&=& 0,\\
\label{come22}
e_{22}(z,\bar{z})&=& i\frac{(\bar{z}-z)}{(\bar{z}+z)},\\
\label{come12}
e_{12}(z,\bar{z})&=& \frac{1}{2}\frac{(z^{2}+\bar{z}^{2})}{(z+\bar{z})^{2}},
\end{eqnarray} 
while the parameter $\Lambda$ from \eqref{airy} is 
\begin{equation}
\label{comLam}
\Lambda(z,\bar{z})=2i\frac{(z^{2}-\bar{z}^{2})}{(z^{2}+\bar{z}^{2})}.
\end{equation}
The differential equation \eqref{airy} after insertion of \eqref{comLam} and rearrangement reduces to
\begin{equation}
\label{airy2}
z^{2}\chi_{,zz}+\bar{z}^{2}\chi_{,\bar{z}\bar{z}}=0.
\end{equation}
Recall that a subscript comma denotes partial differentiation.

A solution to \eqref{airy2} is obtained in the particular class that satisfies
\begin{equation}
\label{airy3}
z^{2}\chi_{,zz}(z,\bar{z})=J(z,\bar{z})=-\bar{z}^{2}\chi_{,\bar{z}\bar{z}},
\end{equation}
where
\begin{equation}
\label{Jdef}
J(z,\bar{z}):= z\bar{z}+k,
\end{equation}
for given  non-negative constant $k$.  Observe  that  
\begin{equation}
\label{J2der}
J_{,zz}(z,\bar{z})=J_{,\bar{z}\bar{z}}(z,\bar{z})=0.
\end{equation} 
\begin{rem}[Generalisation]\label{gensolairy}
The class of solutions satisfying \eqref{airy3} may be extended when the definition of $J(z,\bar{z})$ is generalised to
\begin{equation}
  \label{genJdef}
J(z,\bar{z})= \sum_{j=1}^{m}\frac{(z\bar{z})^{j}}{j!}+k,
\end{equation}
where $m$ is a positive integer and k is a prescribed non-negative constant.  %Further brief details are presented in Appendix~\ref{airyderiv}.
\end{rem}

Integration of \eqref{airy3} with $J(z,\,\bar{z})$ given by \eqref{Jdef} and \eqref{genJdef} is presented in the Appendices. For \eqref{Jdef},   the solution is
\begin{equation}
\label{airy4}
\chi(z,\bar{z})=\left(z\bar{z}-k\right)\log{\frac{z}{\bar{z}}},\qquad z\ne 0,
\end{equation}
%When $z=0$ the function $\chi$ becomes arbitrary and therefore still may  be given by \eqref{airy4}.
which in the previously introduced polar coordinates \textcolor{black}{$(r,\theta)$} becomes 
\begin{equation}
\label{airypol}
\chi(x_{1},\,x_{2})= (r^{2}-k)2i\theta.
\end{equation} 
The corresponding stress and  strain components and Lam\'{e} parameters become
%We have with the previously introduced polar coordinates:
\begin{eqnarray}
\label{str11}
\sigma_{11} &=& -2\left[2\theta +\frac{(r^{2}+k)}{r^{2}}\sin{2\theta}\right],\\
\label{str22}
\sigma_{22}&=& 2\left[-2\theta +\frac{(r^{2}+k)}{r^{2}}\sin{2\theta}\right],\\
\label{str12}
\sigma_{12}&=& 2\frac{(r^{2}+k)}{r^{2}}\cos{2\theta},\\
\label{stra11}
e_{11}&=& 0,\\
\label{stra22}
e_{22}&=& \tan{\theta},\\
\label{stra12}
e_{12}&=& \frac{\cos{2\theta}}{4\cos^{2}{\theta}},\\
\label{muexp}
\mu&=& 4\frac{(r^{2}+k)}{r^{2}}\cos^{2}{\theta},\\
\label{lambexp}
\lambda&=& -4\theta\cot{\theta}-4\frac{(r^{2}+k)}{r^{2}}\cos^{2}{\theta},
\end{eqnarray}
which are singular when $r=0,\,\theta=\pm\pi/2.$

Formulae \eqref{muexp} and \eqref{lambexp}  imply  $\mu >0$ and 
\begin{equation}
\label{lambmugen}
\lambda +2\mu =-4\theta\cot{\theta}+4\frac{(r^{2}+k)}{r^{2}}\cos^{2}{\theta}.
\end{equation}
Choose $k=0$ to give $\mu> 0$ and 
\begin{eqnarray}
\nonumber
  \lambda+2\mu &=& 4\cos{\theta}\left(\cos{\theta}-\frac{\theta}{\sin{\theta}}\right)\\
  \label{neglam}
&<& 0,
\end{eqnarray}
because $\sin{\theta}<\theta$ for $-\pi/2\le \theta <\pi/2$.

Now select $k>0$ and consider a point on the curve $\psi_{c}=1$  for which we have $r=2c\cos{\theta}$.  The relation  \eqref{lambmugen} may be written 
\begin{eqnarray}
\nonumber
\lambda+2\mu &=& 4\cos{\theta}\left(\cos{\theta}-\frac{\theta}{\sin{\theta}}\right)+4\frac{k}{r^{2}}\cos^{2}{\theta}\\
\nonumber
&=& 4\cos{\theta}\left(\cos{\theta}-\frac{\theta}{\sin{\theta}}\right)+\frac{k}{c^{2}}\\
\label{kexpr}
&=& 4\theta\cot{\theta}\left(\frac{\sin{2\theta}}{2\theta}-1\right)+\frac{k}{c^{2}}.
\end{eqnarray}
But $k$ is arbitrary and can be chosen sufficiently large to ensure that 
\begin{equation}
\label{poslam}
\lambda +2\mu >0,
\end{equation}
which is the strong ellipticity condition.

%It is instructive to calculate the corresponding traction acting on the curves $\psi_{c}=1$. The task is  undertaken in Section~\ref{trac}.

\begin{rem}The choice \eqref{Jdef} together with its generalisations provide several examples in nonhomogeneous  isotropic compressible plane  elasticity for which Lam\'{e} parameters are strongly elliptic but for which the displacement %is non-unique due to the 
possesses ambiguous limiting  behaviour at the origin.
\end{rem}

To further examine the problem under consideration, we calculate the total force and couple acting on $\Omega$  together with the strain energy   when $J$  is given by \eqref{Jdef}.

%\subsection{Mechanical point defect}\label{mech}
%It is  unlikely that the point singularity occurring in the stress and strain components \eqref{str11}-\eqref{stra12} corresponds to those  due to a dislocationusually located at an interior point and caused by self-stress. Neither condition applies in the present problem. 

%The next section employs the strain energy to explore several aspects of this feature for the particular displacement boundary value problem considered in this section. 

%%%%%%%%%%%%%%%%%%%%%ok220524

\section{Total force and couple}\label{fse}

The total force and couple  are obtained by the customary procedure of first deriving these quantities for the region $\Omega$ punctured by the deletion  of a neighbourhood of the origin. Contracting to zero the diameter of this neighbourhood yields the desired  results. The additional notation required is now introduced.

\subsection{Geometry}\label{fgeom}

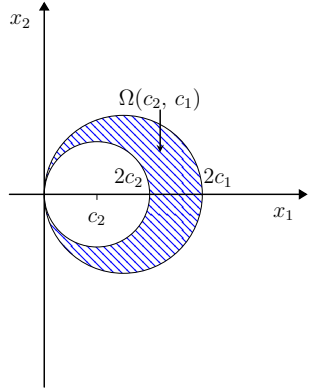
\begin{figure}[htb]
\centerline{
\resizebox{6cm}{!}{
\def\firstcircle{(1,0) circle (1.0cm)}
\def\secondcircle{(1.5,0) circle (1.5cm)}
\begin{tikzpicture}[
  line cap=round,
  line join=round,
  >=Triangle,
  myaxis/.style={->,thick}]

\begin{scope}
        \begin{scope}[even odd rule]% first circle without the second
            \clip \firstcircle (-3,-3) rectangle (3,3);
        \fill[pattern=north west lines,pattern color=blue] \secondcircle;
        \end{scope}
        \draw \firstcircle;
        \draw \secondcircle;
    \end{scope}

\node[above] at (3.3,0) {$2c_1$};

\node[above] at (1.6,0) {$2c_2$};

\node at (2.2,1.8) {$\Omega(c_{2},\,c_{1})$};

\draw [thick,-stealth](2.2,1.6) -- (2.2,0.8);

\draw[myaxis] (-0.66,0) -- (5.00,0) node[below left = 1mm] {$x_1$};      
\draw[myaxis] (0,-3.66) -- (0,3.66) node[below left = 1mm] {$x_2$};

\draw (1,0)--(1,-.1) node[below = 1mm] {$c_2$};
\end{tikzpicture}
}
}
\caption{\label{figomega}The region $\Omega(c_{2},\,c_{1})$.}
\end{figure}

Let $\Omega(c_{2},\,c_{1}),\, 1\le c_{2}<c_{1}\le R,$ be the region bounded internally by the  circle $\psi_{c_{2}}(x_{1},x_{2})=1$ and externally by $\psi_{c_{1}}(x_{1},x_{2})=1$ where $\psi_{c}$ is given by \eqref{phidef}. See Figure \ref{figomega}. Let the bounding circles  be denoted by $\partial\Omega_{c_{\alpha}},\,\alpha=1,2.$  Let $B(0,a),\,a\geq 0$ denote the ball of radius $a$ centred at the origin which is at the common point of the intersection of the family of circles  $\psi_{c}(x_{1},x_{2})=1,\,1\le c\le R.$ That is
\begin{equation}
\label{bdef}
B(0,a):=\left\{(x_{1},\,x_{2}): x^{2}_{1}+x_{2}^{2}\le a^{2}\right\}.
\end{equation}

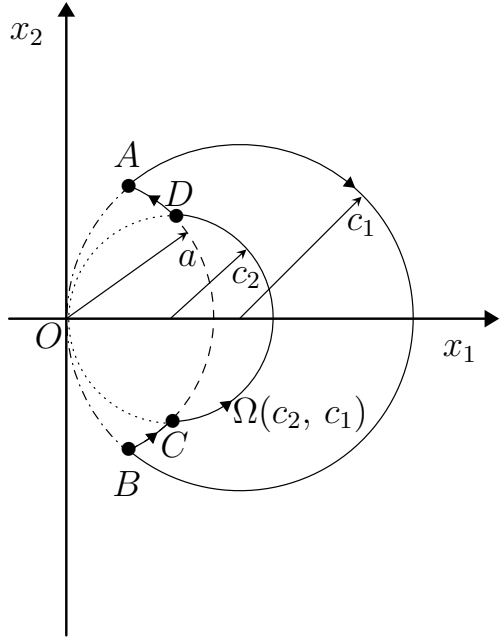
\begin{figure}[htb]
\centerline{
\resizebox{7cm}{!}{
\begin{tikzpicture}[ % scale=5,
  line cap=round,
  line join=round,
  >=Triangle,
  myaxis/.style={->,thick}]

\node at (2.7,-1.1) {$\Omega(c_{2},\,c_{1})$};

\node at (0.72,1.52) {\textbullet};
\node at (1.27,1.17) {\textbullet};
\node at (1.3,1.45) {$D$};
\node at (0.7,1.9) {$A$};
\node at (0.72,-1.52) {\textbullet};
\node at (1.23,-1.19) {\textbullet};
\node at (1.25,-1.45) {$C$};
\node at (0.7,-1.9) {$B$};

\node at (-0.2,-0.2) {$O$};

% arrows decoration
% very primitive way!

% inner circle 

\begin{scope}[ decoration={
  markings,
  mark=at position 0.2 with {\arrow{>}}}
 ] 

\draw[postaction={decorate}] (1.23,-1.19) arc [radius=1.2cm, start angle=272, end
angle= 448];

\end{scope}

% small arcs

\begin{scope}[decoration={
  markings,
  mark=at position 0.7 with {\arrow{>}}}
 ] 

\draw[postaction={decorate}] (0.72,-1.52) arc [radius=1.7cm, start angle=295, end
angle= 314];

\draw[postaction={decorate}] (1.27,1.17) arc [radius=1.7cm, start angle=45, end
angle= 65];

\end{scope}

% outer circle

\begin{scope}[decoration={
  markings,
  mark=at position 0.7 with {\arrow{<}}}
 ]

\draw[postaction={decorate}] (0.72, -1.52) arc [radius=2.0cm, start angle=230, end
angle= 490];

\end{scope}

\draw[dotted] (1.27,1.19) arc [radius=1.2cm, start angle=88, end
angle= 272];

\draw[dash dot]  (0.72, 1.52) arc [radius=2.0cm, start angle=130, end
angle= 230];

\draw[dashed] (0.72,-1.52) arc [radius=1.7cm, start angle=295, end
angle= 405];

% \draw [thick,-stealth](2.2,1.6) -- (2.2,0.8);
\draw [-stealth](1.2,0) -- (2.08,0.8) node[below = 0.6mm] {$c_2$};
\draw [-stealth](2.0,0) -- (3.41,1.41) node[below = 0.6mm] {$c_1$};
\draw [-stealth](0,0) -- (1.4,1.00) node[below = 0.6mm] {$a$};

\draw[myaxis] (-0.66,0) -- (5.00,0) node[below left = 1mm] {$x_1$};      
\draw[myaxis] (0,-3.66) -- (0,3.66) node[below left = 1mm] {$x_2$};

\end{tikzpicture}
}
}
\caption{\label{fig1} $\partial \Omega(c_2,\,c_1) \cap
\partial B(0,a)$.}
\end{figure}

Define the punctured  region $\Omega_{a}(c_{2},\,c_{1})$ by
\begin{equation}
\label{Ompudef}
\Omega_{a}(c_{2},\,c_{1}):=\Omega(c_{2},\,c_{1})\backslash \left\{\Omega(c_{2},\,c_{1})\cap B(0,a)\right\},\qquad a\geq 0,
\end{equation}
%where for completeness we note that
%\begin{equation}
\label{bdef}
%B(0,a):=\left\{(x_{1},x_{2}):x_{1}^{2}+x_{2}^{2} \le a^{2}\right\}.
%\end{equation}
The boundary of the punctured  region  consists of the curves 
\begin{eqnarray}
\label{partbdry}
\partial\Omega_{a}^{c_{\alpha}}&:=& \partial\Omega_{c_{\alpha}}\backslash\left\{\partial\Omega_{c_{\alpha}}\cap B(0,a)\right\}, \qquad \alpha=1,2,\\
\label{Bup}
\partial_{AD}B(0,a)&:=&\left\{(x_{1},x_{2}): (x_{1},x_{2})\in \partial B(0,a)\cap \Omega(c_{2},c_{1}), x_{2}>0\right\},\\
\label{Bdown}
\partial_{CB}B(0,a)&:=&\left\{(x_{1},x_{2}): (x_{1},x_{2})\in \partial B(0,a)\cap \Omega(c_{2},c_{1}), x_{2}<0\right\},
\end{eqnarray}
(Path integrals over $\partial\Omega^{c_{2}}_{a}$ and $\partial_ {AD}B(0,a),\,\partial_{CB}B(0,a)$ will be taken anti-clockwise, while that over $\partial\Omega^{c_{1}}_{a}$ will be taken clockwise.)
%%%%%%%%%%%%%%%%%%%%ok220524
%\subsection{Intersection of $\partial B(0,a)$ and $\partial\Omega(c_{2},c_{1})$.}\label{intersect}
The Cartesian coordinates $(x_{1},\,x_{2})$ of the intersection points $A, D,C, B$ of $\partial B(0,a)$ with $\partial \Omega(c_{2},c_{1})$ (see Figure \ref{fig1}) %have Cartesian coordinates $(x_{1},\,x_{2})$  that a
are  solutions to
\[
x_{1}^{2}+x_{2}^{2}=a^{2},
\]
\[
x_{1}^{2}-2x_{1}c_{\alpha}+x_{2}^{2}=0.
\]

%\begin{figure}[htb]
%\centerline{\resizebox{10cm}{!}{\includegraphics{fig1}}}
%\caption{\label{fig1} $\partial\Omega(c_{2},\,c_{1})\cap\partial B(0,a)$}
%\end{figure}

\begin{table}[h]
\begin{center}
\begin{tabular}{|c|c|c|c|c|c|}\hline
Point &$x_{1}$&$x_{2}$& $\tan{\theta}$& $\cos{\theta}$&$\sin{\theta}$\\ \hline\hline
A & $a^{2}/2c_{1}$ & $(a/2c_{1}) (4c^{2}_{1}-a^{2})^{1/2}$ & $a^{-1}(4c^{2}_{1}-a^{2})^{1/2}$ & $a/2c_{1}$ & $(1/2c_{1}) (4c^{2}_{1}-a^{2})^{1/2}$\\ \hline
D & $a^{2}/2c_{2}$ & $(a/2c_{2})(4c^{2}_{2}-a^{2})^{1/2}$ & $a^{-1}(4c^{2}_{2}-a^{2})^{1/2}$ & $a/2c_{2}$ & $(1/2c_{2}) (4c^{2}_{2}-a^{2})^{1/2}$\\ \hline
B &  $a^{2}/2c_{1}$ & $-(a/2c_{1}) (4c^{2}_{1}-a^{2})^{1/2}$ & $-a^{-1}(4c^{2}_{1}-a^{2})^{1/2}$ & $a/2c_{1}$ & $-(1/2c_{1}) (4c^{2}_{1}-a^{2})^{1/2}$\\ \hline
C & $a^{2}/2c_{2}$ & $-(a/2c_{2}) (4c^{2}_{2}-a^{2})^{1/2}$ & $-a^{-1}(4c^{2}_{2}-a^{2})^{1/2}$ & $a/2c_{2}$ & $-(1/2c_{2}) (4c^{2}_{2}-a^{2})^{1/2}$\\ \hline
\end{tabular}
\caption{Points of intersection}\label{intersect}
\end{center}
\end{table}

 Table \ref{intersect}  lists these  coordinates %and  polar coordinates related by  $x_{1}=r\cos{\theta},\,x_{2}=r\sin{\theta}$, of the points of intersection. Also displayed are
 and also  trigonometric functions of the corresponding polar coordinate $\theta$ where, for example, the polar coordinates of point $A$ are $(r_{A},\,\theta_{A}).$
%%%%%%%%%%%%%%%%%%%ok220524
%\vspace{0.4cm}
%%%%%%%%%%%%%55error220524

%\medskip
%\
%\begin{center}
%Table 1
%\end{center}

%\medskip

Obviously, $a$ must satisfy
\[
a\le 2\le  2c_{2},
\]
whereas
\[
\theta_{A}=-\theta_{B},\qquad \theta_{D}=-\theta_{C}.
\]

Moreover, it is a trivial observation that $a\rightarrow 0$ implies
\[
x_{1}\rightarrow 0,\qquad x_{2}\rightarrow 0,
\]
and
\[
\theta_{A}\rightarrow \pi/2, \qquad \theta_{D}\rightarrow \pi/2,\qquad
\theta_{B}\rightarrow -\pi/2,\qquad \theta_{C}\rightarrow -\pi/2.
\]
%whereas when $c_{1}\rightarrow c_{2}$ we have
%%\begin{eqnarray}
%\theta_{A}&\rightarrow& \theta_{D}\\
%\theta_{B}&\rightarrow&\theta_{C}.
%\end{eqnarray}

%%%%%%%%%%%%%%%%%%%%%%%%%error220524

%%%%%%%%%%%%%%%%%%%%%%%%040624here
\subsection{Traction and couple}\label{trac}
%Section~\ref{asf} establishes that the vector field with components \eqref{u1} and \eqref{u2} is a  displacement vector field in a linear nonhomogeneous isotropic compressible elastic body defined on the lens  shaped region  $\Omega$ and in equilibrium under zero body force and (different) rigid body displacement boundary conditions. The corresponding Lam\'{e} parameters are specified by \eqref{muexp} and \eqref{lambexp}, with the necessarily compatible linear strains expressed by \eqref{stra11}-\eqref{stra12} and stress  by \eqref{str11}-\eqref{str12}. Apart from $e_{11}$ these components all become indeterminate at the origin where the boundary of $\Omega$ is cusp-like. As argued in Section~\ref{base}, this singular behaviour may be due to rupture.

%Most point defects discussed in the literature are located at an interior point and the  elastic body is homogeneous and isotropic.  Moreover, (see \cite{v05}), an interior  point force defect is shown to be equilibrated by tractions over any closed curve or surface surrounding the defect with  the corresponding displacement being continuous. By contrast,  for a dislocation  the displacement exhibits a jump discontinuity while  the resultant traction vanishes over  any surrounding closed curve. In either case, a simple scaling argument shows that the stress and strain components behave like  $O(r^{-1}) $  where $r$ is the distance from the point defect. %as $r\rightarrow 0$ and the point defect is approached. .

We consider the punctured region $\Omega_{a}(c_{2},\,c_{1})$ and compute the resultant traction and couple acting over $\partial_{AD}B(0,a)\cup \partial_{BC}B(0,a)$ 
or, equivalently,  over $\partial \Omega_{a}^{c_{1}}\cup\partial\Omega_{a}^{c_{2}}$.

 The stress components \eqref{str11}-\eqref{str12} combined with the  components  of the unit outward normal on the circle $\psi_{c}(x)=1$  given by %(see \eqref{n1c1} and \eqref{n2c1})
\begin{equation}
\label{n}
n_{1}=(x_{1}-c)/c,\qquad n_{2}= x_{2}/c,
\end{equation}
provide the traction across the circle. % $\psi_{c}(x_{1},x_{2})=1.$ 
Thus, at the point $(x_{1},\,x_{2})$ of the punctured curve, we have
\begin{eqnarray}
\nonumber
F_{1}(x_{1},x_{2})&=& \sigma_{11}n_{1}+\sigma_{12}n_{2}\\
\nonumber
&=& -2\left[2\theta+\left(1+\frac{k}{r^{2}}\right)\sin{2\theta}\right]\frac{(x_{1}-c)}{c}+2\frac{x_{2}}{c}\left(1+\frac{k}{r^{2}}\right)\cos{2\theta}\\
\nonumber
&=&-2\left[2\theta+\left(1+\frac{k}{r^{2}}\right)\sin{2\theta}\right]\cos{2\theta}+2\left(1+\frac{k}{r^{2}}\right)\cos{2\theta}\sin{2\theta}\\
\label{trac1}
&=& -4\theta \cos{2\theta},\\
\nonumber
F_{2}(x_{1},x_{2})&=& \sigma_{21}n_{1}+\sigma_{22}n_{2}\\
\nonumber
&=& 2\left(1+\frac{k}{r^{2}}\right)\cos{2\theta}\frac{(x_{1}-c)}{c}+2\left[-2\theta+\left(1+\frac{k}{r^{2}}\right)\sin{2\theta}\right]\frac{x_{2}}{c}\\
\nonumber
&=& 2\left(1+\frac{k}{r^{2}}\right)\cos^{2}{2\theta} +2\left[-2\theta+\left(1+\frac{k}{r^{2}}\right)\sin{2\theta}\right]\sin{2\theta}\\
\label{trac2}
&=& -4\theta \sin{2\theta}+2\left(1+\frac{k}{r^{2}}\right).
\end{eqnarray}

These expressions indicate the type of singularity occurring at the origin where  $\theta=\pm\pi/2,\, r=0$. On other  parts of the boundary, the traction is non-singular and continuous. Inspection of \eqref{trac1} indicates that  $F_{1}(x_{1},\,x_{2})$ is   an odd function of $\theta$ in the range $-\pi/2\le \theta \le \pi/2$, so that the  total traction acting on $\Omega$ in the $x_{1}$-direction is zero.  
%%%%%%%%%%%%%%%%%%%%%%%%5555Tract

The resultant traction in the $x_{2}$ -direction is obtained from expression  \eqref{trac2} integrated over the boundary curves \eqref{partbdry} with the unit outward normal on $\psi_{c_{1}}=1$  having components 
\begin{eqnarray}
  \label{n1c1}
  n_{1}&=&(x_{1}-c_{1})/c_{1}= -u_{2}/c_{1},\\
  \label{n2c1}
  n_{2}&=&x_{2}/c_{1}=u_{1}/c_{1},
  \end{eqnarray}
while those on $\psi_{c_{2}}=1$ are
\begin{eqnarray}
  \label{n1c2}
  n_{1}&=& -(x_{1}-c_{2})/c_{2}=u_{2}/c_{2},\\
  \label{n2c2}
  n_{2}&=& -x_{2}/c_{2}=-u_{1}/c_{2}.
  \end{eqnarray}

%by \eqref{n1c1}, \eqref{n}.% and \eqref{n1c2},\eqref{n2c2}.

Let us also recall from \eqref{fampol2} and \eqref{phitheta} that in polar coordinates $(r,\,\theta)$  the circle $\psi_{c_{\alpha}}(x_{1},\,x_{2})=1$ is given by
\begin{equation}
  \label{cural}
    r=2c_{\alpha}\cos{\theta},
\end{equation}
and the line element by
\begin{equation}
  \label{line}
ds=c_{\alpha}d(2\theta).
\end{equation}

It follows that  the integral of $F_{2}(.)$ over $\partial\Omega_{a}^{c_{1}}$ can be written 
\begin{eqnarray}
\label{F2c1}
\int_{\partial\Omega_{a}^{c_{1}}}F_{2}(.)\,ds &=& \int^{\theta_{B}}_{\theta_{A}}\left[-4\theta\sin{2\theta}+2\left(1+\frac{k}{r^{2}}\right)\right]c_{1}\,d(2\theta)\\
\nonumber
&=& 2c_{1}\int^{\theta_{B}}_{\theta_{A}}2\theta\, d \left(\cos{2\theta}\right)+4c_{1}\int^{\theta_{B}}_{\theta_{A}}\,d\theta+\frac{k}{c_{1}}\int^{\theta_{B}}_{\theta_{A}} \sec^{2}{\theta}\,d\theta\\
\nonumber
&=& 4c_{1}\left(\theta_{B}\cos{2\theta_{B}}-\theta_{A}\cos{2\theta_{A}}\right)-2c_{1}\left(\sin{2\theta_{B}}-\sin{2\theta_{A}}\right)+4c_{1}\left(\theta_{B}-\theta_{A}\right)\\
\nonumber
& & +\frac{k}{c_{1}}\left(\tan{\theta_{B}}-\tan{\theta_{A}}\right)\\
\label{restract21}
&=& -8c_{1}\theta_{A}\cos{2\theta_{A}}+4c_{1}\sin{2\theta_{A}}-8c_{1}\theta_{A} -2\frac{k}{c_{1}}\tan{\theta_{A}}.
\end{eqnarray}

For the integral over $\partial\Omega_{a}^{c_{2}}$, the unit outward normal is given by \eqref{n1c2} and \eqref{n2c2}. Integration, conducted anti-clockwise, leads to
\begin{equation}
\label{F2c2}
\int_{\partial\Omega_{a}^{c_{2}}}F_{2}(.)\,ds = \int^{\theta_{D}}_{\theta_{C}}\left[-4\theta\sin{2\theta}+2\left(1+\frac{k}{r^{2}}\right)\right]c_{2}\,d(2\theta),
\end{equation}
which  by comparison with \eqref{F2c1} yields
\begin{equation}
\label{restract22}
\int_{\partial\Omega_{a}^{c_{2}}}F_{2}(.)\,ds=8c_{2}\theta_{D}\cos{2\theta_{D}}-4c_{2}\sin{2\theta_{D}}+8c_{2}\theta_{D} +2\frac{k}{c_{2}}\tan{\theta_{D}}.
\end{equation}

When $a>0$, let $T_{2a}$ denote the resultant force in the $x_{2}$-direction acting over the surface $\partial_{DA}B(0,a)\cup\partial _{BC}B(0,a)$. The punctured region is in equilibrium under zero body force so that
\[
T_{2a}+\int_{\partial\Omega_{a}^{c_{1}}}F_{2}(.)\,ds+\int_{\partial\Omega_{a}^{c_{2}}}F_{2}(.)\,ds=0,
\]
and therefore  from \eqref{F2c1} and \eqref{restract22} we obtain  
%Addition of \eqref{restract21} and \eqref{restract22} leads when $a>0$ to  the resultant force $T_{2a}$ in the $x_{2}$-direction acting on the singularity. We have 
\begin{eqnarray}
\nonumber
T_{2a}&= &-8\left(c_{2}\theta_{D}\cos{2\theta_{D}}-c_{1}\theta_{A}\cos{2\theta_{A}}\right)-4\left(c_{1} \sin{2\theta_{A}}-c_{2}\sin{2\theta_{D}}\right)
-8\left(c_{2}\theta_{D}-c_{1}\theta_{A}\right)\\
\label{resF2a}
& & -2k\left(\frac{\tan{\theta_{D}}}{c_{2}}-\frac{\tan{\theta_{A}}}{c_{1}}\right).
\end{eqnarray}

By reference to Table 1 for sufficiently small $a$, we have
\begin{eqnarray}
  \nonumber
  \tan{\theta_{A}}&=& \left(2c_{1}/a\right)\left(1-a^{2}/4c^{2}_{1}\right)^{1/2}\\
  \nonumber
  &=&\left(\frac{2c_{1}}{a}\right)\left[1-\frac{a^{2}}{8c_{1}^{2}}-\frac{1}{8}\left(\frac{a^{2}}{4c_{1}^{2}}\right)^{2}\ldots\right].
\end{eqnarray}
Similarly, we have
\[
\tan{\theta_{D}}=\left(\frac{2c_{2}}{a}\right)\left[1-\frac{a^{2}}{8c_{2}^{2}}-\frac{1}{8}\left(\frac{a^{2}}{4c_{2}^{2}}\right)^{2}\ldots\right],
\]
which leads to the limits
\begin{eqnarray}
  \label{tanlim}
  \lim_{a\rightarrow 0}{|\tan{\theta_{A}}-\left(2c_{1}/a\right)|}&=&0,\\
 \label{tanDlim}
\lim_{a\rightarrow 0}{|\tan{\theta_{D}}-\left(2c_{2}/a\right)|}&=&0.
\end{eqnarray}

Now let $a\rightarrow 0$ so that $\theta_{A}\rightarrow \pi/2,\,\theta_{D}\rightarrow \pi/2$.   In view of \eqref{tanlim} and \eqref{tanDlim}, the resultant traction $T_{2}$ in the $x_{2}$-direction becomes
\begin{equation}
\label{resF2}
T_{2}=0,\qquad a\rightarrow 0.
\end{equation} 

We next consider the anti-clockwise moment, or couple,   of the traction at a point $(x_{1},\,x_{2})$  of $\psi_{c}=1$.  On using \eqref{trac1} and \eqref{trac2}, we have 
\begin{eqnarray}
\nonumber
 \left(x_{1}F_{2}-x_{2}F_{1}\right)
&=& c(\cos{2\theta}+1)\left[-4\theta\sin{2\theta}+2\left(1+\frac{k}{r^{2}}\right)\right]+4c\theta\sin{2\theta}\cos{2\theta}\\
\nonumber
&=& -4c\theta \sin{2\theta}+4c\cos^{2}{\theta}\left(1+\frac{k}{r^{2}}\right)\\
\label{pointmom}
&=& -4c\theta\sin{2\theta}+2c\left(1+\cos{2\theta}\right)+\frac{k}{c}.
\end{eqnarray} 

Thus, the resultant anti-clockwise moment $\Gamma^{c_{1}}_{a}$ acting over the curve  $\partial\Omega^{c_{1}}_{a}$ from \eqref{pointmom} is
\begin{eqnarray}
\nonumber
\Gamma^{c_{1}}_{a}&=& -2c_{1}\int^{\theta_{B}}_{\theta_{A}}2\theta\sin{2\theta}\,c_{1}d(2\theta)+2c_{1}\int^{\theta_{B}}_{\theta_{A}}\left(1+\cos{2\theta}\right)\,c_{1}d(2\theta)+2k(\theta_{B}-\theta_{A})\\
\nonumber
&=& 4c^{2}_{1}(\theta_{B}\cos{2\theta_{B}}-\theta_{a}\cos{2\theta_{A}})+2(k+2c^{2}_{1})(\theta_{B}-\theta_{A})
\\
\label{momtc1}
&=& -8c^{2}_{1}\theta_{A}\cos{2\theta_{A}}-4(k+2c^{2}_{1})\theta_{A},
\end{eqnarray} 
while the clockwise resultant moment over $\partial\Omega_{a}^{c_{2}}$ becomes
\begin{equation}
\label{momtc2}
\Gamma^{c_{2}}_{a}=8c^{2}_{2}\theta_{D}\cos{2\theta_{D}}+4(k+2c^{2}_{2})\theta_{D}.
\end{equation}
Because  the elastic body occupying $\Omega_{a}(c_{2},\,c_{1})$ is in equilibrium, the total moment  $\Gamma_{a}$ satisfies
\begin{equation}
\nonumber
\Gamma_{a}+\Gamma^{c_{1}}_{a}+\Gamma^{c_{2}}_{a}=0,\qquad a>0.
\end{equation}
In consequence from \eqref{momtc1} and \eqref{momtc2}, we obtain
\begin{equation}
\Gamma_{a} =8(c^{2}_{1}\theta_{A}\cos{2\theta_{A}}-c^{2}_{2}\theta_{D}\cos{2\theta_{D}})+4(k+2c^{2}_{1})\theta_{A}-4(k+2c^{2}_{2})\theta_{D}.
\end{equation}
Let $a\rightarrow 0$ and put $c_{2}=1,\,c_{1}=R$, to conclude that  the total moment  $\Gamma$ acting on $\Omega$ is
\begin{equation}
\label{momt}
\Gamma=4(R^{2}-1)\pi.
\end{equation}

The results of this Section, consistent with the interpretation of rupture or tearing due to twisting of the region $\Omega$,  prove the next Proposition:

\begin{proposition}\label{totalload}
  The total force acting   on the region $\Omega$ is zero  but  the resultant couple is of magnitude $4(R^{2}-1)\pi$.
\end{proposition}

The next Section examines the total strain energy.% of the elastic body occupying the region $\Omega$ in equilibrium subject to zero body force and the same rigid body boundary conditions.

\section{The strain energy}\label{gense}

%We further examine the singularity at the origin by calculating the total strain energy 
%of the elastic body occupying the region $\Omega$ subject to the same Dirichlet boundary %conditions.

Define   strain energy in the region $\Omega_{a}(c_{2},c_{1})$  to be 
\begin{equation}
\label{gensepreg}
2E_{a}(c_{2},c_{1}):= \int_{\Omega_{a}(c_{2},c_{1})}\sigma_{\alpha\beta}e_{\alpha\beta}\,dx_1 dx_2,
\end{equation}
where the strain components $e_{\alpha\beta}$ are those in \eqref{stra11}-\eqref{stra12}, and the stress components, obtained from the  Airy stress function \eqref{airypol}, are given in \eqref{str11}-\eqref{str12}. On recalling the stress equilibrium equations \eqref{seq} and  the symmetry of the stress tensor, we   apply the divergence theorem  to obtain
\begin{eqnarray}
\nonumber
2E_{a}(c_{2},c_{1})&=& (1/2)\int_{\Omega_{a}(c_{2},c_{1})}(\sigma_{\alpha\beta}u_{\alpha,\beta}+\sigma_{\alpha\beta}u_{\beta,\alpha})\,dx_{1}dx_{2}\\
\nonumber
&=& \int_{\Omega_{a}(c_{2},c_{1})}\sigma_{\alpha\beta}u_{\alpha,\beta}\,dx_{1}dx_{2}\\
\nonumber
&=& \oint_{\partial\Omega_{a}(c_{2},c_{1})}\sigma_{\alpha\beta}u_{\alpha}n_{\beta}\,ds\\
\label{enboundar}
&=& \oint_{\partial\Omega_{a}(c_{2},c_{1})}\left[(u_{1}\sigma_{11}+u_{2}\sigma_{21})n_{1}+(u_{1}\sigma_{12}+u_{2}\sigma_{22})n_{2}\right]\,ds,
\end{eqnarray}
where the curvilinear integral is taken  clockwise, respectively anticlockwise, as appropriate  and $n_{\alpha}$ are components of the  unit outward normal on $\partial\Omega_{a}(c_{2},c_{1})$.

In general,  the Lam\'{e} parameters  do not explicity appear in  the expression \eqref{enboundar}. 

The aim is  to evaluate the integral \eqref{enboundar}, compute the limit
\[
\lim_{a\rightarrow 0} {E_{a}(c_{2},c_{1})},
\]
and then set $c_{2}=1,\,c_{1}=R$.

We commence by separately computing   \eqref{enboundar}   on the respective  parts of the boundary. This  requires decomposition  into four components:
\begin{equation}
\label{decomint}
\oint_{\partial\Omega_{a}(c_{2},c_{1})}(.)\,ds=\int_{\partial\Omega_{a}^{c_{1}}}(.)\,ds  +\int_{\partial_ {BC}B(0,a)}(.)\,ds +\int_{\partial\Omega^{c_{2}}_{a}}(.)\,ds+\int_{\partial_{DA}B(0,a)}(.)\,ds
\end{equation}
which for convenience is written as
\begin{equation}
\label{vwdef}
\oint_{\partial\Omega_{a}(c_{2},c_{1})}(.)\,ds =V_{1}+W_{1}+V_{2}+W_{2}.
\end{equation}

%Cartesian components of the  unit outward normal on $\psi_{c_{1}}=1$ are
%\begin{eqnarray}
%\label{n1c1}
%n_{1} &=&  (x_{1}-c_{1})/c_{1} = -u_{2}/c_{1},\\
%\label{n2c1}
%n_{2} &=& x_{2}/c_{1} =u_{1}/c_{1},
%\end{eqnarray}
%while for those on $\psi_{c_{2}}=1$, we have
%\begin{eqnarray}
%\label{n1c2}
%n_{1} &= &-(x_{1}-c_{2})/c_{2} =  u_{2}/c_{2},\\
%\label{n2c2}
%n_{2} &=& -x_{2}/c_{2} =- u_{1}/c_{2}.
%\end{eqnarray}
%On $\partial B(0,a)$ the components of the unit outward normal are
%\begin{equation}
%\label{nball}
%n_{1}= -x_{1}/a=-\cos{\theta},\qquad n_{2}=-x_{2}/a=-\sin{\theta}.
%\end{equation}

%In polar coordinates the  curve $\psi_{c_{\alpha}}=1$ is given by
%\[
%r=2c_{\alpha}\cos{\theta},
%\] 
%and  the line element by
%\[
%ds =c_{\alpha}d(2\theta).
%\]

%The terms in the integrand are determined by \eqref{trac1} and \eqref{trac2} when the  stress is derived from the complex formulation for which the strain energy is computed in the next section. 

%%%%%%%%%%%%%%%%%%%%%%%140524

On conversion to polar coordinates $(r,\,\theta)$, the Cartesian  components of the displacement  at the point  $(x_{1},\,x_{2})$ on the curve  $ \psi_{c_{\alpha}}(x_{1},\,x_{2})=1$ are
\begin{eqnarray}
\label{uen1}
u_{1}&=& x_{2}=r\sin{\theta}=c_{\alpha}\sin{2\theta},\\
\label{uen2}
u_{2}&=& (x^{2}_{2}-x^{2}_{1})/2x_{1}=(c_{\alpha}-x_{1})=-c_{\alpha}\cos{2\theta},
\end{eqnarray}
and on $\partial B(0,a)$ are
\begin{eqnarray}
\label{uenb1}
u_{1}&=& a\sin{\theta},\\
\label{uenb2}
u_{2} &=& -a\cos{2\theta}/(2\cos{\theta}),
\end{eqnarray}
where we note that components of the unit outward normal are
\begin{equation}
  \label{nball}
  n_{1}=-x_{1}/a=-\cos{\theta},\qquad n_{2}=-x_{2}/a=-\sin{\theta}.
  \end{equation}

On using \eqref{n1c1}, \eqref{n2c1}, \eqref{str11}-\eqref{str12},  and integrating in the clockwise sense,  we obtain
\begin{eqnarray}
\nonumber
V_{1}&=& \int_{\partial\Omega^{c_{1}}_{a}}u_{\alpha}\sigma_{\alpha\beta}n_{\beta}\,ds\\
\nonumber
&=& \int_{\theta_{A}}^{\theta_{B}}\left[(u_{\alpha}\sigma_{\alpha 1})(-u_{2}/c_{1}) +(u_{\alpha}\sigma_{\alpha 2})(u_{1}/c_{1})\right]c_{1}\,d(2\theta)\\
\nonumber
&=& \int_{\theta_{A}}^{\theta_{B}}\left[u_{1}u_{2}(\sigma_{22}-\sigma_{11})+\sigma_{12}(u^{2}_{1}-u^{2}_{2})\right]\,d(2\theta)\\
\nonumber
&=& 2\int_{\theta_{A}}^{\theta_{B}}\left(1+k/r^{2}\right)\left[2u_{1}u_{2}\sin{2\theta}+(u^{2}_{1}-u^{2}_{2})\cos{2\theta}\right]\,d(2\theta),
\end{eqnarray}
which on substitution from \eqref{uen1} and \eqref{uen2} yields
\begin{eqnarray}
\nonumber
V_{1}&=& 2c^{2}_{1}\int_{\theta_{A}}^{\theta_{B}}\left(1+k/r^{2}\right)\left[-2\sin^{2}{2\theta}\cos{2\theta} +(\sin^{2}{2\theta}-\cos^{2}{2\theta})\cos{2\theta}
\right]\,d(2\theta)\\
%\nonumber
%&=& -2c^{2}_{1}\int_{\theta_{A}}^{\theta_{B}}\left(1+k/r^{2}\right)\cos{2\theta}\,d(2\theta)\%\
%\nonumber
%&=& -2c^{2}_{1}\int^{\theta_{B}}_{\theta_{A}}\cos{2\theta}\,d(2\theta)+k\int_{\theta_{A}}\\\\%^{\theta_{B}}(2-\sec^{2}{\theta})\,d\theta\\
\nonumber
&=& 2c^{2}_{1}(\sin{2\theta_{A}}-\sin{2\theta_{B}})+2k(\theta_{B}-\theta_{A})+k(\tan{\theta_{A}}-\tan{\theta_{B}})\\
\label{enc1}
&=& 4c^{2}_{1}\sin{2\theta_{A}}-4k\theta_{A}+2k\tan{\theta_{A}}.
\end{eqnarray}
The evaluation of $V_2$ is similar. In view of \eqref{n1c2} and \eqref{n2c2}, and on taking  the integral anti-clockwise, we have
\begin{eqnarray}
\nonumber
V_2&=&\int_{\partial\Omega^{c_{2}}_{a}}u_{\gamma}\sigma_{\gamma\delta}n_{\delta}\,ds\\
&=& 2c^{2}_{2}\int_{\theta_{C}}^{\theta_{D}}(1+k/r^{2})\cos{2\theta}\,d(-2\theta)\\
%\nonumber
%&=& 2c^{2}_{2}(\sin{2\theta_{C}}-\sin{2\theta_{D}})\\
%\nonumber
%& & \qquad +2k(\theta_{D}-\theta_{C}) -k(\tan{\theta_{D}}-\tan{\theta_{C}})\\
\label{enc2}
&=& -4c^{2}_{2}\sin{2\theta_{D}}+4k\theta_{D}-2k\tan{\theta_{D}}.
\end{eqnarray}

Addition of \eqref{enc1} and \eqref{enc2} gives  in the notation of \eqref{vwdef}
\begin{equation}
\label{21V}
V_{1}+V_{2}=4c^{2}_{1}\sin{2\theta_{A}}-4c^{2}_{2}\sin{2\theta_{D}}-4k(\theta_{A}-\theta_{D})+2k\left(\tan{\theta_{A}}-\tan{\theta_{D}}\right).
\end{equation}

Before treating  the integrals $W_{1},\,W_{2}$, we examine the limiting behaviour of $V_{1}+V_{2}$ as $a\rightarrow 0.$ We have
\begin{equation}
\nonumber
\lim_{a\rightarrow 0}{\left(V_{1}+V_{2}\right)}= 0-0 +2k\lim_{a\rightarrow 0}{\left(\tan{\theta_{A}}-\tan{\theta_{D}}\right)}.
\end{equation}
%But on reference to Table 1 we further conclude for sufficiently small $a$ that 
%\begin{eqnarray}
%\nonumber
%\tan{\theta_{A}}&=& \left(2c_{1}/a\right)\left(1-a^{2}/4c^{2}_{1}\right)^{1/2}\\
%\nonumber
%&=& \left(2c_{1}/a\right)\left[1-\frac{a^2}{8c^{2}_{1}}-%frac{1}%{8}\left(\frac{a^2}%{4c^{2}%_{1}}\right)^{2}\ldots\right]\\
%\end{eqnarray}
%and 
%\begin{eqnarray}
%\nonumber
%\tan{\theta_{D}}&=& \left(2c_{2}/a\right)\left(1-a^{2}/4c^{2}_{2}\right)^{1/2}\\
%\nonumber
%&=& \left(2c_{2}/a\right)\left[1-\frac{a^{2}}{8c_{2}^{2}}-\frac{1}{8}\left(\frac{a^{2}}{4c%_{2}^{2}}\right)^{2}\dots\right]
%\end{eqnarray}
%so that in the limit we have
%\begin{eqnarray}
%\label{tanlim}
%\lim_{a\rightarrow 0}{|\tan{\theta_{A}}- \left(2c_{1}/a\right)|&=& 0\\
%\label{tanDlim}
%\lim_{a\rightarrow 0}{|\tan{\theta_{D}}-\left(2c_{2}/a\right)|}&=&0.
%\end{eqnarray}
which on recalling \eqref{tanlim} and \eqref{tanDlim}, we may write 
\begin{equation}
\label{Vlimexp}
\lim_{a\rightarrow 0}{|\left(V_{1}+V_{2}\right)|}- (4k/a)|\left(c_{1}-c_{2}\right)|=0.
\end{equation}

The corresponding integrands of $W_{1},\,W_{2}$  over appropriate parts of $\partial B(0,a)$ are evaluated as follows:
\begin{eqnarray}
\nonumber
u_{1}\sigma_{11}+u_{2}\sigma_{12}&=& -2a\sin{\theta}\left[2\theta+\left(1+k/a^{2}\right)\sin{2\theta}\right]-a\cos^{2}{2\theta}\left(1+k/a^{2}\right)\sec{\theta}\\
\nonumber
&=& -4a\theta\sin{\theta}-a\left(1+k/a^{2}\right)\left[2\sin{\theta}\sin{2\theta}+\cos^{2}{2\theta}\sec{\theta}\right]\\
\label{F1b}
&=& -4a\theta\sin{\theta}-a\left(1+k/a^{2}\right)\sec{\theta}.
\end{eqnarray}
\begin{eqnarray}
\nonumber
u_{1}\sigma_{12}+u_{2}\sigma_{22}&=& 2a\left(1+k/a^{2}\right)\sin{\theta}\cos{2\theta}-a\sec{\theta}\cos{2\theta}\left[-2\theta +\left(1+k/a^{2}\right)\sin{2\theta}\right]\\
&=& 2a\theta\sec{\theta}\cos{2\theta}+a\left(1+k/a^{2}\right)\left[2\sin{\theta}\cos{2\theta}-\cos{2\theta}\sin{2\theta}\sec{\theta}\right]\\
\label{F2b}
&=& 2a\theta \sec{\theta}\cos{2\theta}.
\end{eqnarray} 

Consequently,  \eqref{nball} enables us  to obtain for $(x_{1},\,x_{2})\in\partial B(0,a)$ the expressions 
\begin{eqnarray}
\nonumber
(u_{1}\sigma_{11}+u_{2}\sigma_{12})n_{1}+(u_{1}\sigma_{12}+u_{2}\sigma_{22})n_{2}&=& a\cos{\theta}\left[4\theta\sin{\theta}+\left(1+k/a^{2}\right)\sec{\theta}\right]\\
\nonumber
&& -2a\theta\sin{\theta}\sec{\theta}\cos{2\theta}\\
\nonumber
&=& a\left(1+k/a^{2}\right)+2a\theta\sec{\theta}\left[\sin{2\theta}\cos{\theta}-\cos{2\theta}\sin{\theta}\right]\\
\label{integrball}
&=& 2a\theta\tan{\theta}+a\left(1+k/a^{2}\right).
\end{eqnarray}

Let $\overline{\partial B}(0,a):= \Omega(c_{2},\,c_{1})\cap \partial B(0,a),$ and take integration over $\overline{\partial B}(0,a)$ in the anti-clockwise sense. Then
%\begin{equation}
%\label{balint}
%\theta_{A}=-\theta_{B},\qquad \theta_{D}=-\theta_{C}.
%\end{equation}
%Then we have
\begin{eqnarray}
\nonumber
\int_{\overline{\partial B}(0,a)}\sigma_{\alpha\beta}u_{\alpha}n_{\beta}\,ds&=& \int_{\overline{\partial B}(0,a)}\left\{2a\theta\tan{\theta}+a\left(1+k/a^{2}\right)\right\}\,ds\\
\nonumber
&=& \int_{\overline{\partial B}(0,a)}\left\{2a\theta\tan{\theta}+a\left(1+k/a^{2}\right)\right\}\,ad(\theta)\\
\nonumber
&=& 2a^{2}\int_{\overline{\partial B}(0,a)}\theta\tan{\theta}\,d\theta+a^{2}\left(1+k/a^{2}\right)\left(\theta_{A}-\theta_{D}+\theta_{C}-\theta_{B}\right)\\
&=& 2a^{2}\int_{\theta_{C}}^{\theta_{B}}\log{\cos{\theta}}\,d\theta+2a^{2}\int_{\theta_{A}}^{\theta_{D}}\log{\cos{\theta}}\,d\theta\\
%\nonumber
%& & \qquad  +2a^{2}\left[\log{\cos{\theta_{A}}\cos{\theta_{C}}-\log{\cos{\theta_{B}}%\cos{\theta_{D}}\right]\
%\nonumber
& & \qquad -2a^{2}\left( \theta\log{\cos{\theta}}\right)|^{\theta_{A}}_{\theta_{D}}-a^{2}\left( \theta\log{\cos{\theta}}\right)|^{\theta_{C}}_{\theta_{B}}
\\
\label{ballev}
& & \qquad +a^{2}\left(1+k/a^{2}\right)\left(\theta_{A}-\theta_{D}+\theta_{C}-\theta_{B}\right),
\end{eqnarray}
where we have used  an  integration by parts, the relation  
\begin{equation}
\nonumber
\tan{\theta}=-\frac{d}{d\theta}\left(\log\cos{\theta}\right),
\end{equation} 
and where
\begin{equation}
\label{ball}
\int_{\overline{\partial B}(0,a)}\log{\cos{\theta}}\,d\theta=\int^{\theta_{A}}_{\theta_{D}}\log{\cos{\theta}}\,d\theta+\int^{\theta_{C}}_{\theta_{B}}\log{\cos{\theta}}\,d\theta.
\end{equation}
%But  reference to Table 1 shows that 
%\begin{eqnarray}
%\nonumber
%\left[\log{\cos{\theta_{A}}\cos{\theta_{C}}-\log{\cos{\theta_{B}}\cos{\theta_{D}}\right]&=&\log{\frac{\cos{\theta_{A}}\cos{\theta_{C}}}{\cos{\theta_{B}}\cos{\theta_{D}}}\\
%\nonumber
%&=& \log 1\\
%\nonumber
%&=& 0,
%\end{eqnarray}
Consequently \eqref{ballev} reduces to 
\begin{eqnarray}
\nonumber
\int_{\overline{\partial B}(0,a)}\sigma_{\alpha\beta}u_{\alpha}n_{\beta}\,ds&=& -2a^{2}\int_{\overline{\partial B}(0,a)}\log{\cos{\theta}}\,d\theta\\
\nonumber
& & \qquad +4a^{2}\left[\theta_{A}\log{\cos{\theta_{A}}}-\theta_{D}\log{\cos{\theta_{D}}}\right]\\
\label{ballfin}
& & \qquad +2a^{2}\left(1+k/a^{2}\right)\left(\theta_{D}-\theta_{A}\right).
\end{eqnarray}
But  \eqref{ball}, the mean value theorem, and the condition $c_{1}>c_{2}$ imply that
\[
\lim_{a\rightarrow 0} \int_{\overline{\partial B}(0,a)}\log{\cos{\theta}}\,d\theta=0,
\]
and thus we  conclude, by reference to Table 1, that 
\begin{eqnarray}
\nonumber
\lim_{a\rightarrow 0}\int_{\overline{\partial B}(0,a)}\sigma_{\alpha\beta}u_{\alpha}n_{\beta}\,ds
\nonumber
&=& 2k\lim_{a\rightarrow 0}\left(\theta_{D}-\theta_{A}\right)+4\lim_{a\rightarrow 0}{a^{2}\left[\theta_{A}\log{\cos{\theta_{A}}}-\theta_{D}\log{\cos{\theta_{D}}}\right]}
\\
\nonumber
&=& 4\lim_{a\rightarrow 0}{a^{2}\left[\theta_{A}\log{\cos{\theta_{A}}}-\theta_{D}\log{\cos{\theta_{D}}}\right]}\\
&=& 4\lim_{a\rightarrow 0}{a^{2}\left[\theta_{A}\log{(a/2c_{1})}-\theta_{D}\log{(a/2c_{2})}\right]}\\
\nonumber
&=& 4\lim_{a\rightarrow 0}{a^{2}\left[(\theta_{A}-\theta_{D})\log{(a/2c_{1})}+\theta_{D}\log{(c_{2}/c_{1})}\right]}\\
\label{Wlim}
&=& 0,
\end{eqnarray}
since for $\gamma >1$ we have $x^{\gamma}\log{x}\rightarrow 0$ as $x\rightarrow 0$.

Hence, we have
\[
\lim_{a\rightarrow 0}{(W_{1}+W_{2})}=0,
\]
and  on appeal to  \eqref{vwdef} and \eqref{Vlimexp}, we are led to 
\begin{equation}
\label{enk0}
\lim_{a\rightarrow 0}{|E_{a}(c_{2},\,c_{1}) -(2k/a)\left(c_{1}-c_{2}\right)|}=0,
\end{equation}
or
 \begin{equation}
\label{ensing}
 E_{a}(c_{2},\,c_{1})=O(a^{-1}).
\end{equation}

The total strain energy of the whole domain $\Omega$ is obtained either by setting $c_{1}=R,\, c_{2}=1$, or additively,  so that  \eqref{enk0} yields
\begin{equation}
\label{toten}
\lim_{a\rightarrow 0}{|E_{a}(1,R) -\left(2k/a\right)(R-1)|}=0.
\end{equation}

When $k>0$ and is sufficiently large, \eqref{poslam} indicates that  the Lam\'{e} parameters are strongly-elliptic and the strain energy at the origin  behaves like $r^{-1}$. When $k=0$ the Lam\'{e} parameters by \eqref{neglam} cease to be strongly-elliptic but no singularities now occur in the strain energy which indeed vanishes by \eqref{enk0}.

\section{Further remarks}\label{concl}
We comment  on some  generalisations and open problems.
\begin{enumerate}
\item A deeper study of the solution presented in this article may clarify the nature of the singular behaviour. We illustrate this in the simpler problem of the Laplace equation $\Delta v = 0$ in an unbounded plane sector $D = \{(r \cos \theta,r \sin\theta) : \theta \in (0,\alpha), r>0) \}$ of opening angle $\alpha\in (0,2\pi)$, cp.~\cite{g85,k67}. Dirichlet boundary conditions $v(r, 0) = v_0$, $v(r \cos\alpha,r \sin\alpha)=v_1$, for $r>0$, are imposed at $\partial D$, with given constants $v_0$, $v_1$. This problem admits a unique solution of finite Dirichlet energy, given by $v = v_0+(v_1-v_0) \frac{\theta}{\alpha}$. However, further solutions of infinite Dirichlet energy are obtained as $v+V$, $V$ being a linear combination of $\log(r)$ and $r^{\frac{\pi k}{\alpha}} \sin\left(\frac{\pi k \theta}{\alpha}\right)$, $k \in \mathbb{Z}\setminus \{0\}$. Uniqueness of such infinite-energy solutions may be recovered by imposing appropriate growth conditions at $r=0$, respectively $r=\infty$. 

\item Let the plane (hollow) domain $\Omega$ be formed from the family of space filling simple closed curves
\begin{equation}
\label{scurve}
\psi_{c,d}(x_{1},\,x_{2})=1,
\end{equation}
where $c,\,d$ are prescribed positive constants. For instance, let
\begin{equation}
\label{circles}
\psi_{c,d}:= \frac{(x_{1}-c-d)^{2}+x^{2}_{2}}{c^2},
\end{equation}
which represent a family of circles centred at $(c+d,\,0)$ and of radius $c$. The domain  $\Omega$ is bounded externally by the  circle $\psi_{c_{1},0}$ and internally by the circle $\psi_{c_{2},\,d_{2}}$ where $0<c_{2}\le c < c_{1},\,0\le d\le d_{2}.$ When $d_{2}>0$, each member of the family intersects the $x_{1}$-axis at distinct points $x_{1}=d,\,d+2c$, whereas when $d =d_{2}=0$ the family of circles \eqref{circles} reduces to that previously considered in Section~\ref{props} and therefore each member intersects the $x_{1}$-axis at  points  $(0,0)$ and $(2c,0)$.
The unit tangent vector $t$ to a circle belonging to the family \eqref{circles} is given by
\begin{eqnarray}
\label{t1cir}
t_{1}&=& \frac{x_{2}}{c},\\
\label{t2cir}
t_{2}&=& \frac{(c+d-x_{1})}{c}.
\end{eqnarray} 
At the two points of intersection of the circle with the $x_{1}$-axis  the unit tangent has components $(0,\,\pm 1)$. The values are maintained as $d\rightarrow 0$. Consequently, the vector field $(ct_{1},\,ct_{2})$ is uniquely defined at all points of the $x_{1}$-axis provided $d_{2}>0.$ It is only when $d_{2}=0$ and $\Omega$ becomes lens shaped with symmetric cusps at the origin that the vector field fails to be unambigously defined   at the origin. 
\item As proposed by Dafermos \cite{d12}, the vector field $u(x_{1},\,x_{2})$ may be continued into the region $\Omega(0,1)$ upon defining the respective components to be 
  \begin{eqnarray}
    \nonumber
   u_{1}(x_{1},\,x_{2}) &=& x_{2},\qquad (x_{1},\,x_{2})\in \Omega(0,1),\\
    \nonumber
    u_{2}(x_{1},\,x_{2}) &=& (1-x_{1}),\qquad  (x_{1},\,x_{2})\in\Omega(0,1).
 \end{eqnarray}
  The continuation corresponds to a rigid  body rotation of $\Omega(0,1)$.

  Although the vector field similarly may be continued into the external region $\Omega(R,\infty)$, the implications await investigation.

  \item When $c_{2}=1$ for the inner boundary $\partial\Omega_{c_{2}}$ is replaced by $c_{2}=\epsilon$ for small but positive $\epsilon$, the lens shaped region $\Omega(\epsilon,R)$ becomes a region pierced by a small hole whose boundary is rigidly rotated. We  likewise may take $R\rightarrow \infty$ to  obtain apparently analogous problems in the positive half-plane  $x_{1}\ge 0$ with singularities distributed along  the  $x_{1}=0$ axis.   

%\item A  possible link awaits investigation between the family of circles \eqref{circles} and the closed loop of dislocations discussed  by Orowan \cite{o34}; see also Nabarro \cite[p.5]{n67}.
\item The semi-inverse method developed in Section~\ref{asf} may be applied to any prescribed conservative plane vector field $u(x_{1},\,x_{2})$ and the region $\Omega$ defined appropriately.  The Airy stress function continues to satisfy \eqref{airy}, which may be of mixed rather than hyperbolic type, and a solution yields an equilibrium  stress distribution in the absence of body force. 
\item Other solutions exist to the hyperbolic equation \eqref{airy} which lead to different sets of nonhomogeneous Lam\'{e} parameters and corresponding stresses. Nevertheless,  in each case the prescribed displacement has the same behaviour at the origin.

\item
 We comment briefly on another extension of the analysis. Consider the vector field  whose Cartesian components are  tangential to a general family of simple closed curves $\widehat{\psi}_c=1$  for variable constant $c$, and instead of  \eqref{t1} and \eqref{t2} for some constant $\varpi$  put
\begin{eqnarray}
\label{u1gen}
u_{1}&=& \varpi\frac{\partial\widehat{\psi}_{c}}{\partial x_{2}},\\
\label{u2gen}
u_{2}&=&-\varpi\frac{\partial\widehat{\psi}_{c}}{\partial x_{1}}.
\end{eqnarray}
Observe that the vector with components \eqref{u1gen} and \eqref{u2gen} for smooth $\widehat{\psi}_{c}$ is \textcolor{black}{isochoric} when the point $(x_{1},\,x_{2})$ is confined to a given curve of the family. The representation  is  a special case of the   plane  Helmholtz decomposition 
\[
u= \nabla\Psi+\nabla \times \psi, 
\]
for scalar potential function $\Psi$ and vector potential function $\psi$. The component $\nabla\Psi$ is parallel to the normal to level surfaces of $\Psi$.  The expressions \eqref{u1gen} and \eqref{u2gen} are recovered on setting $\Psi=0$ and $\psi=(0,\,0,\,\varpi\widehat{\psi}_{c})$.
%The component $\nabla \Psi$ may be regarded as parallel to the normal to level surfaces of $\Psi$. The relationship of the previous analysis to the Helmholtz decomposition awaits comprehensive investigation.

The vector whose components are \eqref{u1gen} and \eqref{u2gen} for suitable choices of $\varpi$ and $\widehat{\psi}_{c}$ correspond to the displacement vector for point defects in plane elastostatic equilibrium problems. For example, when $\varpi=2$ and $\widehat{\psi}_{1}(x) =\log{(x_{1}^{2}+x_{2}^2)}$. Possible relationships are then with Volterra and other dislocations in homogeneous linear isotropic elasticity.  See \cite{l27}.
\item A different generalisation is derived on setting   $u(x_{1}\,x_{2})=cw(x_{1},\,x_{2})t(x_{1},\,x_{2})$ (see \cite{a83}), where $w(x)$ is an arbitrary  scalar function. Thus, 
\begin{eqnarray}
\label{u11}
u_{1}(x_{1},\,x_{2})&=& w(x)x_{2},\\
\label{u22}
u_{2}(x_{1},x_{2})&=& w(x) (c-x_{1}).
\end{eqnarray}
Nevertheless, for simplicity we have set $w=1$ in the preceding discussion.

\item It is of interest to extend the procedure adopted here to three-dimensions using the Maxwell-Morera functions to generalise the Airy stress function.   
\end{enumerate}

\paragraph{Acknowledgement.}
The authors are grateful for support from the ICMS Research in Groups scheme.

 \begin{appendix}
\section{Airy function. Derivation}\label{airyderiv}

%{\color{red}``The following might be included after comments on differentiation,See text of 22.0424''
This appendix provides the detailed elementary derivation of the solution \eqref{airy4} to equation \eqref{airy3}, subject to \eqref{Jdef} or more generally to \eqref{J2der}, which on integration shows that 
\begin{equation}
\label{genJ}
J(z,\bar{z})= zg(\bar{z})+f(\bar{z}),
\end{equation}
where the arbitrary functions $g(.),\,f(.)$ are chosen \textcolor{black}{below. Other choices are presented in Appendix \ref{Jgengen}}.

%Appendix~\ref{Jgengen}  enlarges the   choice of the function \eqref{genJ}. 

 Computations are given only for the left side of \eqref{airy3}. The right side may be treated similarly. 

Consequently, consider
\[
z^{2}\chi_{,zz}=J(z,\bar{z}),
\]                                                  
which leads to
\begin{eqnarray}
\nonumber
\chi_{,zz}(z,\bar{z})&=& \frac{J(z,\bar{z})}{z^{2}},\qquad z \ne 0\\
\nonumber
&=& -J\frac{d}{dz}\left(\frac{1}{z}\right)\\
\nonumber
&=& -\frac{d}{dz}\left(\frac{J}{z}\right)+\frac{J_{,z}}{z}\\
\nonumber
&=& -\frac{d}{dz}\left(\frac{J}{z}\right)+J_{,z}\frac{d}{dz}\left(\log{z}\right)\\
\label{chi1}
&=& -\frac{d}{dz}\left(\frac{J}{z}\right)+\frac{d}{dz}\left(J_{,z}\log{z}\right),
\end{eqnarray}
on recalling \eqref{J2der}. Integration of the last expression yields 
\[
\chi_{,z}=-\frac{J}{z}+J_{,z}\log{z}+G(\bar{z}),
\] 
where $G(\bar{z})$ is an arbitrary function of $\bar{z}$. Write the last equation  as
\begin{eqnarray}
\nonumber
\chi_{,z}&=& -J\frac{d}{dz}\left(\log{z}\right)+J_{,z}\log{z}+G(\bar{z})\\
\nonumber
&=& -\frac{d}{dz}\left(J\log{z}\right)+2J_{,z}\frac{d}{dz}\left(z\log{z}-z\right)+G(\bar{z})\\
\nonumber
&=&-\frac{d}{dz}\left(J\log{z}\right)+2\frac{d}{dz}\left(J_{,z}\left(z\log{z}-z\right)\right)+G(\bar{z}),
\end{eqnarray}
which on integration gives
\begin{equation}
\label{chi2}
\chi(z,\bar{z})= -J\log{z}+2J_{,z}\left(z\log{z}-z\right)+zG(\bar{z})+F(\bar{z})
\end{equation}
where $F(\bar{z})$ is an arbitrary function of $\bar{z}$.

Substitution for $J$ from \eqref{genJ} in \eqref{chi2} after  simplification gives
\[
\chi(z,\bar{z})= g(\bar{z})z\log{z}-f(\bar{z})\log{z}-2zg(\bar{z})+zG(\bar{z})+F(\bar{z}).
\]
Set
\begin{equation}
g(\bar{z})=k_{1}\bar{z},\qquad f(\bar{z})=k_{2},\qquad G(\bar{z})=2g(\bar{z})-k_{1}\bar{z}\log{\bar{z}},\qquad F(\bar{z})= k_{2}\log{\bar{z}},
\end{equation}
which on putting $k_{1}=1,\,k_{2}=k$ leads to \eqref{airy4} as required.

\section{Generalised $J(z,\bar{z})$}\label{Jgengen}
Let  $J$ be given by
\begin{equation}
\label{Jext}
J(z,\bar{z})=\sum_{j=1}^{m}\frac{(z\bar{z})^{j}}{j!}+k,\qquad m=1,2,3,\ldots
\end{equation}
from which follows
\begin{equation}
\label{Jextde}
\frac{d^{(m+1)}J(z,\bar{z})}{dz^{(m+1)}}=\frac{d^{(m+1)}J(z,\bar{z})}{d\bar{z}^{(m+1)}}=0.
\end{equation}
The procedure outlined in Appendix~\ref{airyderiv} may now be repeated on noting the formulae
\begin{equation}
\label{genparts}
\frac{d}{dz}\left(nz^{n}\log{z}-z^{n}\right)=n^{2}z^{(n-1)}\log{z},\qquad n=1,2,3,\ldots.
\end{equation}

  The case $n=1$ is treated in Appendix~\ref{airyderiv}. The general form of the Airy function $\chi$ may  easily be computed  but the result is not recorded.
\end{appendix}

\end{document}